\documentclass[11pt,a4paper]{amsart}
\usepackage[all]{xy}

\newcommand*{\myref}[1]{(\ref{#1})}
\newtheorem{theorem}{Theorem}[section]
\newtheorem{corollary}[theorem]{Corollary}
\newtheorem{lemma}[theorem]{Lemma}
\newtheorem{proposition}[theorem]{Proposition}

\theoremstyle{definition}
\newtheorem{definition}[theorem]{Definition}
\newtheorem{remark}[theorem]{Remark}
\newcommand{\ie}{{\em i.e.}\ }
\newcommand{\cf}{{\em cf.}\ }
\newcommand{\eg}{{\em e.g.}\ }
\newcommand{\ko}{\: , \;}
\newcommand{\vir}{\:, \;}

\newcommand*{\ul}[1]{\underline{#1}}
\newcommand*{\ol}[1]{\overline{#1}}
\renewcommand*{\tilde}[1]{\widetilde{#1}}
\newdir{ >}{{}*!/-8pt/\dir{>}}
\newdir{ (}{{}*!/-6pt/\dir^{(}}
\newcommand{\Ac}{A^\wedge}

\newcommand{\Aic}{{A_i}^\wedge}

\newcommand{\Ft}{{\tilde{F}}}
\newcommand{\Nl}{{N_\lambda}}
\newcommand{\Qr}{{Q_\rho}}

\newcommand{\Z}{\mathbb{Z}}

\newcommand{\ca}{{\mathcal A}}
\newcommand{\cb}{{\mathcal B}}
\newcommand{\cc}{{\mathcal C}}
\newcommand{\cd}{{\mathcal D}}
\newcommand{\ce}{{\mathcal E}}
\newcommand{\cF}{{\mathcal F}}
\newcommand{\cg}{{\mathcal G}}
\newcommand{\ch}{{\mathcal H}}

\newcommand{\cl}{{\mathcal L}}
\newcommand{\cm}{{\mathcal M}}
\newcommand{\cn}{{\mathcal N}}

\newcommand{\cs}{{\mathcal S}}
\newcommand{\ct}{{\mathcal T}}
\newcommand{\cu}{{\mathcal U}}
\newcommand{\cv}{{\mathcal V}}

\newcommand{\eps}{\varepsilon}
\newcommand*{\opname}[1]{\operatorname{#1}}

\newcommand{\la}{\leftarrow}
\newcommand{\lra}{\longrightarrow}
\newcommand{\lla}{\longleftarrow}
\newcommand{\eqiso}{\stackrel{_\sim}{=}}
\newcommand{\iso}{\stackrel{_\sim}{\rightarrow}}
\newcommand{\inviso}{\stackrel{_\sim}{\leftarrow}}
\newcommand{\Iso}{\stackrel{_\sim}{\longrightarrow}}
\newcommand{\Inviso}{\stackrel{_\sim}{\longleftarrow}}
\newcommand{\id}{\mathbf{1}}

\newcommand{\op}{^{op}}
\newcommand{\ten}{\otimes}

\newcommand{\DG}{\opname{DG}}

\newcommand{\Hom}{\opname{Hom}}

\newcommand{\sk}{\opname{sk}}
\newcommand{\cone}{\opname{cone}}

\newcommand{\Mod}{\opname{Mod}\nolimits}

\newcommand{\coh}{\opname{coh}}

\newcommand{\Add}{\opname{Add}}

\newcommand{\Presh}{\opname{Presh}}
\newcommand{\Sh}{\opname{Sh}}
\newcommand{\colim}{\opname{colim}}
\newcommand{\col}{\opname{col}}

\newcommand{\cok}{\opname{cok}}

\newcommand{\coprodh}{\coprod^{H^0}}
\newcommand{\ta}{\ct^\alpha}
\newcommand{\tb}{\ct^\beta}
\newcommand{\td}{\ct^\delta}

\newcommand{\tg}{\ct^\gamma}

\newcommand{\tl}{\ct^\lambda}

\newcommand{\ts}{\ct^\sigma}

\newcommand{\na}{\cn^\alpha}
\newcommand{\nb}{\cn^\beta}
\newcommand{\nd}{\cn^\delta}
\newcommand{\ngamma}{\cn^\gamma}

\newcommand{\nl}{\cn^\lambda}

\newcommand{\Tagen}{{\langle\ta\rangle}}
\newcommand{\Sgen}{{\langle\cs\rangle}}
\newcommand{\Ggen}{{\langle\cg\rangle}}

\newcommand{\Sa}{{\langle\cs\rangle_\alpha}}
\newcommand{\Sb}{{\langle\cs\rangle_\beta}}
\newcommand{\Sd}{{\langle\cs\rangle_\delta}}
\newcommand{\Sg}{{\langle\cs\rangle}_\gamma}

\newcommand{\Gb}{{\langle\cg\rangle_\beta}}

\newcommand{\QGa}{{\langle Q\cg \rangle_\alpha}}
\newcommand{\QGb}{{\langle Q\cg \rangle_\beta}}

\newcommand{\CA}{{\cc\ca}}

\newcommand{\DA}{{\cd\ca}}
\newcommand{\DAN}{{\cd\ca/\cn}}
\newcommand{\DaA}{{\cd_\alpha\ca}}

\title[Popescu-Gabriel for triangulated categories]{The
  Popescu-Gabriel theorem for triangulated categories}

\author{Marco Porta}

\address{UFR de Math{\'e}matiques, UMR 7586 du CNRS, Case 7012,
   Universit{\'e} Paris 7, 2 place Jussieu, 75251 Paris Cedex 05,
   France}
\email{porta@math.jussieu.fr}

\address{Dipartimento di Matematica ``F. Enriques'', Universit{\`a}
degli Studi di Milano, v. C. Saldini 50, 20133 Milano, Italy}
\email{porta@mat.unimi.it}

\begin{document}

\begin{abstract}
The Popescu-Gabriel theorem states that each Grothendieck abelian
category is a localization of a module category. In this paper, we
prove an analogue where Grothendieck abelian categories are replaced
by triangulated categories which are well generated (in the sense of
Neeman) and algebraic (in the sense of Keller). The role of module
categories is played by derived categories of small differential
graded categories. An analogous result for topological triangulated
categories has recently been obtained by A. Heider.
\end{abstract}

%\begin{classification}
%Primary 18E30; Secondary 16D90.
%\end{classification}

 \date{February 15, 2008}

%\begin{keywords}
%Homological algebra, derived category, homotopy category,
%derived functor, $K$-theory, Hochschild cohomology, Morita theory,
%non commutative algebraic geometry.
%\end{keywords}

\maketitle

\tableofcontents

\section{Introduction}
\label{s:intro}
One of the aims of the present article is to try and answer the
question: what is the analogue, in the realm of triangulated
categories, of the notion of Grothendieck category in the realm of
abelian categories? The best way to proceed seemed to us that of
lifting perhaps the most important theorem involving these notions
from the abelian to the triangulated world. The theorem we are
speaking of is due to Popescu-Gabriel:

\begin{theorem}[Popescu-Gabriel \cite{PopescuGabriel64}]
\label{thm:PG}
Let $\ct$ be a Grothendieck category. Then the following statements are
equivalent:
\begin{itemize}
\item[($i$)] $G\in\ct$ is a generator of $\ct$;
\item[($ii$)] the functor $\Hom(G,-) : \ct \longrightarrow \Mod(A)$,
      where $A=\Hom(G,G)$, is a localization.
\end{itemize}
\end{theorem}

We refer to \cite{Popescu73}, \cite{Takeuchi71}, \cite{Mitchell81},
\cite{Lowen04} for complete proofs of the theorem.

In his book \cite[Def.~1.15, p.~15]{Neeman99}, A. Neeman defined the class of
{\em well generated triangulated categories}. It turns out that this class is a
very good generalization to higher cardinals of the concept of {\em compactly generated triangulated category}. In fact, it preserves the most interesting properties, \eg the validity of the Brown representability theorem \cite{Brown62} and of the Thomason localization theorem \cite[Key Proposition 5.2.2, p.~338]{ThomasonTrobaugh90}, and at the same time
introduces new good features, such as the stability of the new class
under localizations (assuming the quite weak hypothesis that the
kernel of the localization functor is generated by a set of
objects). H. Krause characterized the class of categories introduced
by Neeman as follows \cite{Krause01}. Let $\ct$ be a triangulated
category with suspension functor $\Sigma$ admitting arbitrary
set-indexed coproducts. $\ct$ is well generated in the sense of Krause
\cite{Krause01} if and only if there exists a set $\cg_0$ of objects
with $\Sigma\cg_0=\cg_0$ satisfying the conditions:

\begin{itemize}
\item[(G1)] an object $X\in\ct$ is zero provided that $\ct(G,X)=0$ for
  all $G$ in $\cg_0$;
\item[(G2)] for each family of morphisms $f_i: X_i \to Y_i$, $i\in I$,
  the induced map
\[
\ct(G,\coprod_{i\in I} X_i) \to \ct(G,\coprod_{i\in I} Y_i)
\] is surjective for all $G\in\cg_0$ provided that the maps
\[
\ct(G,X_i) \to \ct(G,Y_i)
\]
are surjective for all $i\in I$ and all $G\in\cg_0$;
\item[(G3)] there is some regular cardinal $\alpha$ such that the
      objects $G\in\cg_0$ are $\alpha${\em -small}, \ie for each
      family of objects $X_i$, $i\in I$, of $\ct$, each morphism
\[
G \to \coprod_{i\in I} X_i
\]
factors through a subsum $\coprod_{i\in J} X_i$ for some subset
$J$ of $I$ of cardinality strictly smaller than $\alpha$.
\end{itemize}

In the case $\alpha=\aleph_0$, the $\aleph_0${\em -compact} objects
are the {\em compact} objects of the classical literature
\cite{Neeman92}, \cite{Neeman96} and the definition of well generated
category reduces to that of {\em compactly generated} one. Well
generated triangulated categories arise very naturally when one
localizes compactly generated ones, as it will be shown in detail in
section~\ref{s:tsubcatloc}. For example, the unbounded derived
category $\cd(\Sh(X))$ of sheaves of abelian groups over a topological
space $X$ is well generated since it is a localization of the derived
category of presheaves $\cd(\Presh(X))$, which is compactly
generated. However, Neeman shows in \cite{Neeman01} that not all
derived categories of sheaves are compactly generated. An example is
the category $\cd(\Sh(X))$ where $X$ is a connected, non compact real
manifold of dimension at least one; in this case, there do not exist
non zero compact objects. In the same article, it is shown that the
derived categories of Grothendieck categories are always well
generated. Another large class of examples arises when one localizes
the derived category $\DA$ of a small DG category $\ca$ at the
localizing subcategory generated by a set of objects. Indeed, since
$\DA$ is a compactly generated triangulated category, such a
localization is always well generated. Now we can state the main
result of this paper. It also gives a positive answer to Drinfeld's
question \cite{Neeman07} whether all well generated categories arise
as localizations of module categories over DG categories, for the class of {\em
algebraic} triangulated categories. Here algebraic means triangle
equivalent to the stable category of a Frobenius category. One can
show that each algebraic triangulated category is triangle equivalent
to a full triangulated subcategory of the category up to homotopy of
complexes over some additive category.

\begin{theorem}
Let $\ct$ be an algebraic triangulated category. Then the following
statements are equivalent:
\begin{itemize}
\item[($i$)] $\ct$ is well generated;
\item[($ii$)] there is a small $\DG$ category $\ca$ such that $\ct$ is
      triangle equivalent to a localization of $\DA$ with respect to a
      localizing subcategory generated by a {\em set} of objects.
\end{itemize}
Moreover, if (i) holds and $\cg \subseteq \ct$ is a full triangulated
subcategory stable under coproducts of strictly less than $\alpha$
factors and satisfying (G1), (G2) and (G3) for some regular cardinal
$\alpha$, the functor
\[
\ct \longrightarrow \Mod\cg \vir X \longmapsto
\Hom_\ct(-,X)|_\cg
\]
lifts to a localization $\ct \longrightarrow \cd(\tilde\cg)$, where
$\tilde\cg$ is a small DG category such that $H^0(\tilde\cg)$ is
equivalent to $\cg$.
\end{theorem}

If $\ct$ is compactly generated, the theorem yields a triangle
equivalence $\ct \longrightarrow \DA$, and we recover Theorem 4.3 of
\cite{Keller94}. Note the structural similarity with the abelian
case. One notable difference is that in the abelian case, one can work
with a single generator whereas in the triangulated case, in general,
it seems unavoidable to use a (small but usually infinite)
triangulated subcategory.

An analogous result for topological triangulated categories has
recently been proved by A. Heider \cite{HeiderThesis}.

\subsection{Organization of the paper}
In section~\ref{s:wgtc}, we present some auxiliary results about well
generated triangulated categories. After recalling the definition
given by Krause (subsection~\ref{ss:defKN}), we establish a small set
of conditions which allows us to show that two well generated
triangulated categories are triangle equivalent
(subsection~\ref{ss:eqwgtc}). 

In section~\ref{s:tsubcatloc}, we recall some basic results about
localizations of triangulated categories and about their thick
subcategories. In subsection~\ref{ss:locwgtc}, we state a theorem
concerning particular localizations of well generated triangulated
categories, those which are triangle quotients by a subcategory {\em
  generated by a set}.

Section~\ref{s:acdc} is the heart of the paper. We construct the
$\alpha${\em-continuous derived category} $\DaA$ of a
{\em homotopically $\alpha$-cocomplete} small DG category $\ca$
(section~\ref{s:acdc}). This construction enjoys a useful and
beautiful property which is the key technical result for proving the
main theorem of the paper: Given a homotopically $\alpha$-cocomplete
pretriangulated DG category $\ca$, we show that its
$\alpha$-continuous derived category $\DaA$ is $\alpha$-compactly
generated by the images of the free DG modules. The proof heavily
uses theorem~\ref{thm:locwgtc} of subsection~\ref{ss:locwgtc} about
localizations of well generated triangulated categories.

The categories $\DaA$ turn out to be the prototypes of the
$\alpha$-compactly generated algebraic DG categories. This
characterization is what we have called the Popescu-Gabriel theorem
for triangulated categories. This is the main result of the paper.
We present it in section~\ref{s:PGforTC}. As an application, we also
give a result about compactifying subcategories of an algebraic well
generated triangulated category. The notion of compactifying
subcategory generalizes that of compactifying generator introduced
by Lowen-Van den Bergh \cite[Ch.~5]{LowenVandenBergh04} in the case
of a Grothendieck abelian category.

\subsection{Acknowledgements} I am very happy to thank Bernhard Keller, my Ph.D. thesis director, for a lot of discussions and helpful comments on this work, as well as for his natural kindness. I also thank Bert van Geemen, the joint tutor of this thesis in Milan, for his interest and support and Amnon Neeman for some comments on a preliminary version of this article. This is the occasion to thank Georges Maltsiniotis for inviting me to present the results of this paper in the working seminar ``Alg{\`e}bre et Topologie Homotopiques'', held in Paris at the ``Institut de Math{\'e}matiques de Jussieu''. I am grateful
to the ``Universit{\'e} Franco-Italienne'', which supported me with a grant.

%\newpage

\section{Well generated triangulated categories}
\label{s:wgtc}

\subsection{Definitions of Krause and Neeman}
\label{ss:defKN}
The notion of well generated triangulated category is due to A. Neeman
\cite[Def.~1.15, p.~15]{Neeman99}. Instead of his original definition, we will use a characterisation due to H. Krause \cite{Krause01} which is closer in
spirit to the definition of Grothendieck abelian categories. We recall
that a {\em regular} cardinal $\alpha$ is a cardinal which is {\em
  not} the sum of fewer than $\alpha$ cardinals, all smaller than
$\alpha$ (see any standard reference about set theory for definitions
and properties of ordinals and cardinals, a very readable one is
\cite{Krivine98}). In this article, we will usually assume that the
cardinals we use are infinite and regular.

\begin{definition}
\label{def:wgtcKrause}
Let $\ct$ be a triangulated category with arbitrary coproducts and
suspension functor $\Sigma$. Let $\alpha$ be an infinite regular
cardinal. Then the category $\ct$ is $\alpha${\em -compactly
  generated} if there exists a set of {\em $\alpha$-good generators}, \ie a set of objects $\cg_0$ such that
$\Sigma\cg_0=\cg_0$, satisfying the conditions:
\begin{itemize}
\item[(G1)] an object $X\in\ct$ is zero if $\ct(G,X)=0$ for all $G$ in
      $\cg_0$;
\item[(G2)] for each family of morphisms $f_i: X_i \to Y_i$, $i\in I$,
  the induced map
\[
\ct(G,\coprod_{i\in I} X_i) \to \ct(G,\coprod_{i\in I} Y_i)
\] is surjective for all $G\in\cg_0$ if the maps
\[
\ct(G,X_i) \to \ct(G,Y_i)
\]
are surjective for all $i\in I$ and all $G\in\cg_0$;
\item[(G3)] all the objects $G\in\cg_0$ are $\alpha${\em -small}, \ie
      for each family of objects $X_i$, $i\in I$, of $\ct$, each
      morphism
\[
G \to \coprod_{i\in I} X_i
\]
factors through a subsum $\coprod_{i\in J} X_i$ for some subset
$J$ of $I$ of cardinality strictly smaller than $\alpha$.
\end{itemize}

A triangulated category is {\em well generated} \cite{Krause01}, if there exists a regular cardinal $\delta$ such that
it is $\delta$-compactly generated.
\end{definition}

Let $\ct$ be a triangulated category with arbitrary coproducts. We
will say that condition (G4) holds for a class of objects $\cg$ of
$\ct$ if the following holds:

\begin{itemize}
\item[(G4)] for each family of objects $X_i$, $i\in I$, of
 $\ct$, and each object $G \in \cg$, each morphism
\[
G \to \coprod_{i\in I} X_i
\]
factors through a morphism $\coprod_{i\in I} \phi_i$: $\coprod_{i \in
  I} G_i \to \coprod_{i \in I} X_i$, with $G_i$ in $\cg$ for all $i
  \in I$.
\end{itemize}

Clearly, condition (G4) holds for the empty class and, if it holds for
a family of classes, then it holds for their union. Thus, for a given
regular cardinal $\alpha$, there exists a unique maximal class
satisfying (G4) and formed by $\alpha$-small objects. Following Krause
\cite{Krause01}, we denote this class, and the triangulated
subcategory on its objects, by $\ta$. Its objects are called the
$\alpha${\em -compact} objects of $\ct$.

\begin{remark}
This definition of $\ta$ is not identical to the one of Neeman
\cite[Def.~1.15, p.~15]{Neeman99}. However, as shown in \cite[Lemma 6]{Krause01}, the two definitions are equivalent if the isomorphism classes of $\ta$ form a set. This always holds when $\ct$ is well generated, \cf \cite{Krause01}.
\end{remark}

In the case $\alpha=\aleph_0$, the $\aleph_0${\em -compact} objects
are the objects usually called compact (also called small). We recall
that an object $K$ of $\ct$ is called {\em compact} if the following
isomorphism holds
\[
\bigoplus_{i \in I}\ct(K,X_i) \iso \ct(K,\coprod_{i \in I}X_i) \ko
\]
where the objects $X_i$ lie in $\ct$ for all $i \in I$, and $I$ is an
arbitrary set. The triangulated category with coproducts $\ct$ is
usually called {\em compactly generated} if condition (G1) holds for a
set $\cg_0$ contained in the subcategory of compact objects $\ct^c =
\ct^{\aleph_0}$. In the case $\alpha=\aleph_0$, the definition of well
generated category specializes to that of compactly generated
category.

Let $\ct$ be a triangulated category with arbitrary coproducts and
$\cg_0$ a small full subcategory of $\ct$. Let $\cg = \Add(\cg_0)$ be
the closure of $\cg_0$ under arbitrary coproducts and direct
factors. A functor $F : \cg\op \to \ca b$ is {\em coherent}
\cite{Cartan53}, \cite{Auslander66} if it admits a presentation
\[
\cg(-,G_1) \to \cg(-,G_0) \to F \to 0
\]
for some objects $G_0$ and $G_1$ of $\cg$. Let $\coh(\cg)$ be the {\em
category of coherent functors} on $\cg$. It is a full subcategory of
the {\em category} $\Mod\cg$ of all additive functors $F : \cg\op \to
\ca b$. Part c) of the following lemma appears in \cite[Lemma
3]{Krause02}, in a version with countable coproducts instead of
arbitrary coproducts. We give a new, more direct proof.

\begin{lemma}
\label{lm:comm}
\
\begin{itemize}
\item[a)] For each object $X$ of $\ct$, the functor $h(X)$ obtained by
restricting $\ct(-,X)$ to $\cg$ is coherent.
\item[b)] The functor $\cg \to \coh(\cg)$ taking $G$ to $h(G)$
      commutes with arbitrary coproducts.
\item[c)] Condition (G2) holds for $\cg_0$ iff $h : \ct \to \coh(\cg)$
  commutes with arbitrary coproducts.
\end{itemize}
\end{lemma}

\begin{proof}
a) We have to show that, for each $X \in \ct$, the functor
$\ct(-,X)|_\cg$, which, {\em a priori}, is in $\Mod\cg$, is in fact
coherent. We choose a morphism $\coprod_{i \in I} G_i \to X \vir G_i
\in \cg_0$, such that each $G \to X \vir G \in \cg_0$, factors through
a morphism $G_i \to X$. Then
\[
\ct(-,\coprod_{i \in I} G_i)|_\cg \to \ct(-,X)|_\cg
\]
is an epimorphism in $\Mod\cg$. We form a distinguished triangle
\[
X' \to \coprod_{i \in I} G_i \to X \to \Sigma X'
\]
in $\ct$. We can continue the construction and choose a morphism
$\coprod_{i \in I'} G'_i \to X' \vir G'_i \in \cg_0$, such that each
$G \to X' \vir G \in \cg_0$, factors through a morphism $G'_i \to
X'$. Then the sequence
\[
\ct(-,\coprod_{i \in I'} G'_i)|_\cg \to \ct(-,\coprod_{i \in I}
G_i)|_\cg \to \ct(-,X)|_\cg \to 0
\]
is a presentation of $\ct(-,X)|_\cg$.

\smallskip
b) Let $(G_i)_{i \in I}$ be a family of objects of $\cg$. We have to
show that the canonical morphism
\[
\coh(\cg)(h(\coprod_{i \in I} G_i), F) \to \prod_{i \in I}
\coh(\cg)(h(G_i), F)
\]
is invertible for each coherent functor $F$. Since $h(G)$ is projective for
each $G$ in $\cg$, it is enough to check this for representable functors $F$.
For these, it follows from Yoneda's lemma and the definition of $\coprod_{i \in
I} G_i$.

\smallskip
c) We suppose that (G2) holds for $\cg_0$.

\smallskip
{\em First step}. For each family $(X_i)_{i \in I}$ of $\ct$, the
canonical morphism $\coprod_{i \in I} h(X_i)
\stackrel{_\varphi}{\longrightarrow} h(\coprod_{i \in I} X_i)$ is an
epimorphism. Indeed, for each $i \in I$, let $G_i \to X_i$ be a
morphism such that
\[
h(G_i) \to h(X_i)
\]
is an epimorphism, where $G_i$ belongs to $\cg$. By b), the functor $h
: \cg \to \coh(\cg)$ commutes with coproducts. Thus, we obtain a
commutative square
\[
\xymatrix{
\coprod_{i \in I} h(G_i) \ar[r] \ar[d]^\wr & \coprod_{i \in I} h(X_i)
\ar[d]^\varphi  \\
h(\coprod_{i \in I} G_i) \ar[r]_<>(0.5)\pi & h(\coprod_{i \in I} X_i).
}
\]
By condition (G2), $\pi$ is an epimorphism. Thus, $\varphi$ is an
epimorphism.

\smallskip
{\em Second step}. For each family $(X_i)_{i \in I}$ of $\ct$, the
canonical morphism $\coprod_{i \in I} h(X_i)
\stackrel{_\varphi}{\longrightarrow} h(\coprod_{i \in I} X_i)$ is an
isomorphism. Indeed, for each $i \in I$, we choose distinguished
triangles
\[
X'_i \to G_i \to X_i \to \Sigma X'_i \ko
\]
and morphisms $G'_i \to X'_i$, where $G_i \to X_i$ is as in the first
step and $G'_i$ belongs to $\cg$, such that
\[
h(G'_i) \to h(X'_i)
\]
is an epimorphism. Then the sequence
\[
0 \to h(\coprod_{i \in I} X'_i) \underset{\iota}{\rightarrow}
h(\coprod_{i \in I} G_i) \underset{\pi}{\rightarrow} h(\coprod_{i \in
  I} X_i) \to 0
\]
is exact. Indeed, coproducts preserve distinguished triangles and $h$
is cohomological since it is the composition of the Yoneda functor
with the restriction functor $F \mapsto F|_\cg$, which is clearly
exact. In particular, $\iota$ is a monomorphism. Since the coproduct
functor $\coprod_{i \in I}$ is right exact, the top morphism of the
square
\[
\xymatrix{
\coprod_{i \in I} h(G'_i) \ar@{->>}[r] \ar[d]^\wr & \coprod_{i \in I}
h(X'_i) \ar@{->>}[d]^\varphi  \\
h(\coprod_{i \in I} G'_i) \ar[r] & h(\coprod_{i \in I}
X'_i)
}
\]
is an epimorphism. By the first step, it follows that the morphism
$\varphi$ is an epimorphism. By b), the morphism $\coprod_{i \in I}
h(G'_i) \to h(\coprod_{i \in I} G'_i)$ is an isomorphism. Therefore,
the morphism
\[
h(\coprod_{i \in I} G'_i) \to h(\coprod_{i \in I} X'_i)
\]
is an epimorphism and the sequence
\[
h(\coprod_{i \in I} G'_i) \to h(\coprod_{i \in I} G_i)
\underset{\pi}{\rightarrow} h(\coprod_{i \in I} X_i) \to 0
\]
is exact. The claim now follows from b) since, $\coprod_{i \in I}$
being a right exact functor, we have a diagram with exact rows
\[
\xymatrix{
\coprod_{i \in I} h(G'_i) \ar[r] \ar[d]^\wr & \coprod_{i \in I} h(G_i)
\ar[r] \ar[d]^\wr & \coprod_{i \in I} h(X_i) \ar[r] \ar[d]^\varphi & 0
\\
h(\coprod_{i \in I} G'_i) \ar[r] & h(\coprod_{i \in I} G_i) \ar[r] &
h(\coprod_{i \in I} X_i) \ar[r] & 0.
}
\]
%\
We suppose now that $h$ commutes with coproducts. We will show that
condition (G2) holds for $\cg$. Let $(f_i : X_i \to Y_i)_{i \in I}$ be
a family of morphisms in $\ct$ such that $\ct(G,f_i) : \ct(G,X_i) \to
\ct(G,Y_i)$ is surjective for all $i \in I$ and all $G \in
\cg_0$. Then $\ct(\coprod_{l \in L} G_l,f_i) : \ct(\coprod_{l \in L}
G_l,X_i) \to \ct(\coprod_{l \in L} G_l,Y_i)$ is surjective for all the
families $(G_l)_{l \in L}$ of $\cg_0$ and all $i \in I$, thanks to the
isomorphisms $\ct(\coprod_{l \in L} G_l,X_i) \to \prod_{l \in L}
\ct(G_l,X_i)$. Moreover, it is trivial to verify that $\ct(A,f_i) :
\ct(A,X_i) \to \ct(A,Y_i)$ is surjective, for all $i \in I$, for each
direct factor $A$ of any object $G \in \cg_0$. Therefore, $\ct(G,f_i)$
is surjective for all $i \in I$ and all $G \in \cg$. Thus,
$\ct(-,X_i)|_\cg \to \ct(-,Y_i)|_\cg$ is an epimorphism for all ${i
  \in I}$. The coproduct $\coprod_{i \in I} \ct(-,X_i)|_\cg \to
\coprod_{i \in I} \ct(-,Y_i)|_\cg$ is still an epimorphism. Since $h$
commutes with coproducts, it follows that $\ct(G,\coprod_{i \in I}
X_i) \to \ct(G,\coprod_{i \in I} Y_i)$ is surjective for all $G \in
\cg$, in particular for all $G \in \cg_0$.
\end{proof}

Consider a triangulated category $\ct$ and a class of its objects
$\cg_0$, satisfying some or all the conditions of
definition~\ref{def:wgtcKrause}. It will be important for us to know
if these conditions continue to hold for different closures of
$\cg_0$.

\begin{proposition}
\label{prop:dense}
Let $\ct$ be a {\em cocomplete} triangulated category, \ie $\ct$
admits all small coproducts. Let $\cg_0$ be a class of objects in
$\ct$, stable under $\Sigma$ and $\Sigma^{-1}$, satisfying conditions
(G2) and (G3) of the definition~\ref{def:wgtcKrause}. Let $\alpha$ be
an infinite cardinal. Let $\cg$ be the closure of $\cg_0$ under
$\Sigma$ and $\Sigma^{-1}$, extensions and $\alpha$-small
coproducts. Then, conditions (G3) and (G4) hold for $\cg$.
\end{proposition}

\begin{proof}
(G3) We directly show that condition (G3) holds for shifts,
$\alpha$-coproducts and extensions of objects in $\cg_0$.
%
% OMITTED PROOF OF THE FOLLOWING ASSERTION
%
% Given an arbitrary object $G$ in $\cg_0$ and a coproduct $\coprod_{i
% \in I}X_i$ of objects in $\ct$, the following sequence of
% isomorphisms, with $n \in \Z$, proves that (G3) holds for shifts:
%
% \begin{eqnarray*}
% \Hom_\ct(\Sigma^n G, \coprod_{i \in I}X_i) & = &
% \Hom_\ct(G,\Sigma^{-n}\coprod_{i \in I}X_i) \\
% & = & \Hom_\ct(G,\coprod_{i \in I}\Sigma^{-n}X_i) \\
% \colim_{J\subset I}\Hom_\ct(\Sigma^n G,\coprod_{i \in J}X_i) & = &
% \colim_{J\subset I}\Hom_\ct(G,\coprod_{i \in J}\Sigma^{-n}X_i) \ko
% \end{eqnarray*}
%
% where $J$ are sets of cardinality strictly smaller than $\alpha$. We
% use the fact that the restriction of $\Sigma$ to $\cg$ is an
% equivalence endofunctor of $\cg$, since this is a triangulated
% subcategory of $\ct$. For the third isomorphism, we use the
% $\alpha$-smallness of $G$.

Since the functor $\Sigma : \ct \to \ct$ is an equivalence, an object
$X$ of $\ct$ is $\alpha$-small iff $\Sigma X$ is $\alpha$-small. Thus,
condition (G3) holds for all objects $\Sigma^n G, G \in \cg_0, n \in
\Z$.

Since $\alpha$-small coproducts commute with $\alpha$-filtered
colimits, condition (G3) holds for $\alpha$-small coproducts of
objects of $\cg_0$. Indeed, let $(G_j)_{j \in J}$, $|J| < \alpha$, be
a family of $\alpha$-small objects of $\cg_0$ and let $(X_i)_{i \in
  I}$ be an arbitrary family of objects of $\ct$. We have the
following sequence of isomorphisms:
\begin{eqnarray*}
\Hom_\ct(\coprod_{j \in J}G_j,\coprod_{i \in I}X_i) & = &
\prod_{j \in J}\Hom_\ct(G_j,\coprod_{i \in I}X_i) \\
& = & \prod_{j \in J}\colim_{I' \subset I}\Hom_\ct(G_j,\coprod_{i \in
  I'}X_i)  \\
\colim_{I' \subset I}\Hom_\ct(\coprod_{j \in J}G_j,\coprod_{i \in
  I'}X_i) & = & 
\colim_{I' \subset I}\prod_{j \in J}\Hom_\ct(G_j,\coprod_{i \in I'}X_i) \ko
\end{eqnarray*}
where the cardinality of the subset $I'$ is strictly smaller than
$\alpha$. The only non trivial isomorphism is the vertical
third which holds since the cardinal $\alpha$ is supposed regular,
hence the colimit is taken over an $\alpha$-filtered set $I$.

Let us consider the (mapping) cone of an arbitrary morphism $G \to G'$
of $\cg_0$
\[
\xymatrix{
G \ar[r] & G' \ar[r] & C \ar[r] & \Sigma G .
}
\]
We can form two long exact sequences by applying the cohomological
functors
\[
\Hom_\ct(-,\coprod_{i \in I}X_i) \quad \mbox{and} \quad \Hom_\ct(-,\coprod_{i
  \in J}X_i)
\]
to the last distinguished triangle. Now we consider the
colimit over the subsets $J \subset I$ of cardinality strictly smaller
than $\alpha$ of the long exact sequence induced by
$\Hom_\ct(-,\coprod_{i \in J}X_i)$. We obtain a long sequence which is
still exact since we are using filtered colimits. There is a natural
map of the two long exact sequences just formed. Let us represent a
part of it in the following diagram, where we write $\col_J$ for
$\colim_{J \subseteq I}$
\[
\xymatrix{
\col_J\ct(G,\coprod_J X_i) \ar[d]^\wr &
\col_J\ct(G',\coprod_J X_i) \ar[d]^\wr \ar[l] &
\col_J\ct(C,\coprod_J X_i) \ar[d] \ar[l] &
\col_J\ct(\Sigma G,\coprod_J X_i) \ar[d]^\wr \ar[l]  \\
\ct(G,\coprod_I X_i) & \ct(G',\coprod_I X_i) \ar[l] &
\ct(C,\coprod_I X_i) \ar[l] & \ct(\Sigma G,\coprod_I X_i)
\ar[l] .
}
\]
The vertical arrows are isomorphisms since $G$, $G'$ are in $\cg_0$,
and we have seen that $\Sigma G$ is $\alpha$-small. Thus, the third
vertical arrow is an isomorphism by the Five-Lemma and $C$ is
$\alpha$-small, too.

\smallskip
(G4) We call $\cu$ the full subcategory of $\ct$ formed by the objects
$X \in \cg$ which satisfy the following condition. Given a morphism
\[
f : X \xymatrix{ \ar[r] &} \coprod_{i\in I}Y_i \ko
\]
where $(Y_i)_{i\in I}$ is a family of objects in $\ct$, there exists a
family $(X_i)_{i\in I}$ of objects of $\cg$ and some morphisms
$\varphi_i : X_i \to Y_i$ such that $f$ factors as in the diagram
%\[
%\xymatrix{
%X \ar[rr]^f \ar[dr] && \coprod_{i\in I} Y_i  \\
%& \coprod_{i\in I} X_i. \ar[ur]_{\coprod_{i\in I}\varphi_i}
%}
%\]
\[
\xymatrix@ur{
X \ar[ddrr]^f \ar[dd]  \\\\
\coprod_{i\in I}X_i. \ar[rr]_{\coprod_{i\in I}\varphi_i} &&
\coprod_{i\in I} Y_i
}
\]
We shall show:
\begin{itemize}
\item[a)] the subcategory $\cu$ contains $\cg_0$;
\item[b)] the subcategory $\cu$ is stable under formation of
      $\alpha$-coproducts;
\item[c)] the subcategory $\cu$ is closed under $\Sigma$, $\Sigma^{-1}$
  and under extensions.
\end{itemize}
It follows by the properties a), b), c) that $\cu = \cg$, which shows
that the condition (G4) holds for $\cg$.

\smallskip
a) Let $G_0$ be an object in $\cg_0$ and $f : G_0 \to \coprod_{i\in
  I}Y_i$ a morphism in $\ct$, where $(Y_i)_{i\in I}$ is a family of
objects in $\ct$. For every $i \in I$, let $(G_{ij} \to Y_i)_{j\in
  J_i}$, where $G_{ij} \in \cg$, be a family of morphisms such that
every morphism $G_0 \to Y_i$ factors through one of the morphisms
$G_{ij} \to Y_i$. Then, the morphism
\[
\varphi_i : \coprod_{j\in J_i}G_{ij} \xymatrix{ \ar[r] &} Y_i
\]
induces a surjection
\[
\Hom_\ct(G_0,\coprod_{j\in J_i}G_{ij}) \xymatrix{ \ar[r] &}
\Hom_\ct(G_0,Y_i) \ko
\]
for every $i \in I$. By (G2), the map
\[
\Hom(G_0,\coprod_{i\in I}\coprod_{j\in J_i}G_{ij}) \xymatrix{ \ar[r]
  &} \Hom(G_0,\coprod_{i\in I}Y_i)
\]
is a surjection. Therefore, there exists a morphism
\[
\tilde f : G_0 \xymatrix{ \ar[r] &} \coprod_{i\in I}\coprod_{j\in
  J_i}G_{ij}
\]
such that the composition
\[
G_0 \xymatrix{ \ar[rr]^{\tilde f} &&} \coprod_{i\in I}\coprod_{j \in
  J_i}G_{ij} \xymatrix{ \ar[rr]^{\coprod_{i \in I}\varphi_i} &&}
\coprod_{i \in I}Y_i
\]
is equal to $f$. We have supposed that $G$ is $\alpha$-small
(condition (G3) holds for $\cg_0$). Therefore the morphism
\[
G_0 \xymatrix{ \ar[r]^{\tilde f} &} \coprod_{i\in I}\coprod_{j\in
  J_i}G_{ij} = \coprod_{(i,j)\in\cl}G_{ij} \ko
\]
where $\cl$ is the set of pairs $(i,j)$ with $i \in I$ and $j \in J_i$, factors through the sub-sum
\[
\coprod_{\stackrel{}{(i,j)\in\Lambda}}G_{ij} =
\coprod_{j\in \tilde I}\coprod_{j\in {\tilde J}_i}G_{ij} \ko
\]
where $\Lambda \subseteq \cl$ is a subset of cardinality strictly smaller than $\alpha$. Let $\tilde I$ be the set of indices $i \in I$ such that $\Lambda$ contains a pair of the form $(i,j)$. Then $\tilde I$ is of cardinality strictly smaller than $\alpha$. Now for each $i \in \tilde I$, let $\tilde{J_i}$ be the set of indices $j \in J_i$ such that $\Lambda$ contains the pair $(i,j)$. Then each $\tilde{J_i}$ is of cardinality strictly smaller than $\alpha$. Now for $i \notin \tilde I$, put ${\tilde J}_i = \emptyset$. Then we have
\[
\coprod_{\stackrel{}{(i,j)\in\Lambda}}G_{ij} =
\coprod_{\stackrel{}{i\in I}}\coprod_{j\in {\tilde J}_i}G_{ij} .
\]
Let $Y_i = \coprod_{j\in {\tilde J}_i}G_{ij}$. Then $f$ factors as
\[
G_0 \xymatrix{ \ar[rr] &&} \coprod_{i\in I}Y_i \xymatrix{
  \ar[rr]^{\coprod_{i\in I}(\varphi_i\mid_{Y_i})} &&} \coprod_{i\in
  I}X_i .
\]
As $|\tilde J|<\alpha$, $Y_i$ lies in $\cg$ for all $i \in I$.

\smallskip
b) Let $(U_j)_{j\in J}$ be a family of $\cu$ where $|J|<\alpha$. Let
\[
f : \coprod_{j\in J}U_j \xymatrix{ \ar[r] &} \coprod_{i\in I}X_i \ko
\]
be a morphism in $\ct$, where $(X_i)_{i\in I}$ is a family of
$\ct$. Let
\[
f_j : U_j \xymatrix{ \ar[r] &} \coprod_{i\in I}X_i
\]
be the component of $f$ associated to $j \in J$. For each $j \in J$,
since $U_j$ lies in $\cu$, there exists a factorization
\[
\xymatrix@ur{
U_j \ar[ddrr]^{f_j} \ar[dd]  \\\\
\coprod_{i\in I}Y_{ji} \ko \ar[rr]_{\coprod_{i\in I}\varphi_{ji}} &&
\coprod_{i\in I} X_i
}
\]
where $(Y_{ji})_{i\in I}$ is a family of $\cg$. Then, we have the
factorization
\[
\coprod_{j \in J}U_j \xymatrix{ \ar[r] &} \coprod_{j \in J}\coprod_{i
  \in I} Y_{ji} \xymatrix{ \ar[r]^{\varphi} &} \coprod_{i \in I}X_i \ko
\]
which we can write as
\[
\coprod_{j\in J}U_j \xymatrix{ \ar[rr] &&} \coprod_{i \in I}\coprod_{j
  \in J} Y_{ji} \xymatrix{ \ar[rr]^{\coprod_{i \in I}\varphi_i} &&}
\coprod_{i\in I}X_i \ko
\]
where $\coprod_{j\in J} Y_{ji}$ belongs to $\cg$ since
$|J|<\alpha$. Therefore, $\coprod_{j\in J}U_j$ lies in $\cu$.

\smallskip
c) Clearly, $\cu$ is stable under the action of $\Sigma$ and
$\Sigma^{-1}$. Let
\[
\xymatrix{
X \ar[r] & X' \ar[r] & X'' \ar[r]  & \Sigma X   \\
}
\]
be a distinguished triangle of $\ct$ such that $X$, $X'$ are in
$\cu$. Let
\[
X'' \xymatrix{ \ar[r]^{f''} &} \coprod_{i\in I}Y_i
\]
be a morphism of $\ct$ where $(Y_i)_{i\in I}$ is a family of $\ct$. We
have the factorization
\[
\xymatrix{
X' \ar[rr] \ar[d]^{f'} && X'' \ar[d]^{f''}   \\
\coprod_{i\in I}X'_i \ar[rr]_{\coprod_{i\in I}\varphi_i} &&
\coprod_{i\in I}Y_i \ko
}
\]
where $X' \in \cg$ and $\varphi_i : X'_i \to Y_i$ are morphisms in
$\ct$ for all $i \in I$. We can extend each $\varphi_i$ to a
distinguished triangle, take coproducts over $I$ and then complete the
square above to a morphism of distinguished triangles (using axiom
TR3 of triangulated categories):
\[
\xymatrix{
X \ar[rr] \ar@{-->}[d]^f && X' \ar[rr] \ar[d]^{f'} && X'' \ar[rr]
\ar[d]^{f''} && \Sigma X \ar@{-->}[d]^{\Sigma f}  \\
\coprod_{i\in I}Z_i \ar[rr] && \coprod_{i\in I}X'_i
\ar[rr]_{\coprod_{i\in I}\varphi_i} && \coprod_{i\in I}Y_i
\ar[rr]_{\coprod_{i\in I}\eps_i} && \Sigma\coprod_{i\in I}Z_i .
}
\]
The objects $X$, $X'$ and $\Sigma X$ belong to $\cu$. Thus, the
morphisms $f$, $f'$ and $\Sigma f$ above factor through a coproduct
taken over $I$ of objects in $\cg$. We have the commutative diagram
\[
\xymatrix{
X \ar[rr] \ar[d] && X' \ar[rr] \ar[d] && X'' \ar[rr] && \Sigma X
\ar[d]  \\
\coprod_{i\in I}X_i \ar[rr]_{\coprod_{i\in I} u_i} \ar[d] &&
\coprod_{i\in I}X'_i \ar@{=}[d] &&&& \Sigma\coprod_{i\in I}X_i \ar[d]
\\
\coprod_{i\in I}Z_i \ar[rr] && \coprod_{i\in I}X'_i
\ar[rr]_{\coprod_{i\in I}\varphi_i} && \coprod_{i\in I}Y_i
\ar[rr]_{\coprod_{i\in I}\eps_i} && \Sigma\coprod_{i\in I}Z_i \ko
}
\]
where the morphisms $u_i : X_i \to X'_i$, $i \in I$, are in
$\cg$. Now, we extend the morphisms $u_i$ to distinguished triangles
and then form the distinguished triangle of coproducts over
$I$. Successively, we form morphisms of distinguished triangles using
axiom TR3 of triangulated categories, obtaining the maps $g$ and
$\coprod_{i\in I}\psi_i$ as shown in the diagram
\[
\xymatrix{
X \ar[rr] \ar[d] && X' \ar[rr] \ar[d] && X'' \ar[rr] \ar@{-->}[d]^g
\ar@/_2.7pc/[dd]_<>(0.3){f''} \ar@/_1.4pc/@{.>}[lldd]_h && \Sigma X
\ar[d] \\
\coprod_{i\in I}X_i \ar[rr]_{\coprod_{i\in I} u_i} \ar[d] &&
\coprod_{i\in I}X'_i \ar[rr] \ar@{=}[d] && \coprod_{i\in I}X''_i
\ar[rr] \ar@{-->}[d]^{\coprod_{i\in I}\psi_i} && \Sigma \coprod_{i\in
  I}X_i \ar[d]   \\
\coprod_{i\in I}Z_i \ar[rr] && \coprod_{i\in I}X'_i
\ar[rr]_{\coprod_{i\in I}\varphi_i} && \coprod_{i\in I}Y_i
\ar[rr]_{\coprod_{i\in I}\eps_i} && \Sigma\coprod_{i\in I}Z_i .
}
\]
Note that the subdiagram between $X''$ and $\coprod_{i\in I}Y_i$ does
{\em not} commute, \ie the composition $(\coprod_I\psi_i) \circ g$ is
in general {\em not} equal to $f''$. Anyway, by composing with
$\coprod_I\eps_i$, we obtain
\[
(\coprod_I\eps_i) \circ (\coprod_I\psi_i) \circ g = (\coprod_I\eps_i)
\circ f'' .
\]
Therefore, by applying $\Hom_\ct(X'',-)$ to the distinguished triangle
in the third row of the last diagram, it is immediate that
\[
(\coprod_I\psi_i) \circ g - f'' = (\coprod_I\varphi_i) \circ h \ko
\]
for some morphism $h : X'' \to \coprod_{i\in I}X'_i$, as in the
diagram. Then, the correct expression of $f''$ is
\[
f'' =  (\coprod_I\psi_i) \circ g + (\coprod_I\varphi_i) \circ (-h)
\]
which shows that $f''$ factors as
\[
X'' \xymatrix{ \ar[rrr]_<>(0.5){\begin{bmatrix} g \\ -h \end{bmatrix}}
  &&&}
%{\left[\stackrel{g}{-h} \right]}
\coprod_{i\in I}X''_i \oplus \coprod_{i\in I}X'_i
\xymatrix{ \ar[rrr]_<>(0.5){[\coprod_I\psi_i \vir \coprod_I\varphi_i]}
  &&&} \coprod_{i\in I}Y_i .
\]
The previous factorization of $f''$ is trivially equivalent to the
following
\[
X'' \xymatrix{ \ar[rrr] &&&} \coprod_{i\in I}(X''_i \oplus X'_i)
\xymatrix{ \ar[rrr]_<>(0.5){\coprod_I[\psi_i \vir \varphi_i]} &&&}
\coprod_{i\in I}Y_i .
\]
Now, $X''_i \oplus X'_i$ is in $\cg$ for all $i\in I$ by
construction.
\end{proof}

There are two immediate and useful corollaries.

\begin{corollary}
\label{coro:thick}
Let $\ct$ be a cocomplete triangulated category. Let $\cg_0$ be a
class of objects in $\ct$ satisfying all the conditions of the last
proposition. Let $\alpha$ be an infinite cardinal. Let $\cg$ be the
closure of $\cg_0$ under $\Sigma$ and $\Sigma^{-1}$, extensions,
$\alpha$-small coproducts and direct factors, \ie $\cg = \langle \cg_0
\rangle_\alpha$ in the notation of~\ref{ss:subcat} below. Then,
conditions (G3) and (G4) hold for $\cg$.
\end{corollary}

\begin{proof}
The proof of the preceding proposition works for (G4) if we verify
that the subcategory $\cu$ is also closed under direct factors, \ie
that it is thick \myref{ss:subcat}.

Let $U$ be an object in $\cu$ and $U = U' \oplus U''$. Then, there is
a section $i$ of the projection $p : U \to U'$. Let $f : U' \to
\coprod_{i\in I}W_i$ be a morphism in $\ct$. The composition $f \circ
p$ factors as
\[
\xymatrix{
U \ar@<+0.35ex>[r]^p \ar[dr]_g & U' \ar@<+0.35ex>[l]^i \ar@{-->}[d]^{g
  \circ i} \ar[r]^<>(0.5)f & \coprod_{i\in I}W_i  \\
& \coprod_{i\in I}V_i \ar[ur]_{\coprod_{i\in I}\phi_i} \ko
}
\]
where the objects $V_i$ are in $\cg$ and the morphisms $\phi : V_i \to
W_i$ in $\ct$, for all $i \in I$. Then $f$ also factors over
$\coprod_{i\in I}V_i$, through the morphism $g \circ i$. Indeed, $f
\circ p = (\coprod_I\phi_i) \circ g$, and $f \circ p \circ i =
(\coprod_I\phi_i) \circ g \circ i$, but $p \circ i$ is the identity
morphism of $U'$.

The proof of the preceding proposition works for (G3) if we verify
that the direct factors of the objects in $\cg_0$ are $\alpha$-small,
too. This requires the construction of a diagram structurally
identical to the one above. Therefore, we omit it.
\end{proof}

\begin{corollary}
\label{coro:thickwgtc}
Let $\alpha$ be an infinite regular cardinal. Let $\ct$ be a
triangulated category $\alpha$-compactly generated by a set
$\cg_0$. Let $\cg$ be the closure of $\cg_0$ under $\Sigma$ and
$\Sigma^{-1}$, extensions and $\alpha$-small coproducts. Let $\langle
\cg_0 \rangle_\alpha$ be the closure of $\cg$ under direct
factors. Then, $\ct$ is $\alpha$-compactly generated by both $\cg$ and
$\langle \cg_0 \rangle_\alpha$.
\end{corollary}

\begin{proof}
The condition (G1) clearly holds for both $\cg$ and $\langle \cg_0
\rangle_\alpha$, since they contain $\cg_0$. The conditions (G3) and
(G4) hold for $\cg$ by proposition~\ref{prop:dense} and for $\langle
\cg_0 \rangle_\alpha$ by corollary~\ref{coro:thick}. Moreover,
condition (G4) easily implies (G2).
\end{proof}

\subsection{Equivalences of well generated triangulated categories}
\label{ss:eqwgtc}
This subsection is devoted to establishing a small set of conditions
which allows us to show that two well generated triangulated
categories are triangle equivalent.

\begin{proposition}
\label{prop:eqwgtc}
Let $\ct$ and $\ct'$ be two triangulated categories admitting arbitrary set-indexed coproducts. Let $\alpha$ be a regular cardinal and $\cg \subset \ct$ and $\cg' \subset \ct'$ two {\em $\alpha$-localizing} subcategories, \ie thick and closed under formation of $\alpha${\em -small} coproducts \myref{ss:subcat}. Suppose that $\cg$ and $\cg'$ satisfy conditions (G1), (G2), (G3) for the cardinal $\alpha$. Let $F : \ct \to \ct'$ be a triangle functor which commutes with all coproducts and induces an equivalence $\cg \to \cg'$. Then $F$ is an equivalence of triangulated categories.
\end{proposition}

\begin{proof} $1^{st}$ {\em step}: The functor $F$ induces an equivalence
\[
\Add\cg \lra \Add\cg' .
\]

As $F$ commutes with coproducts and induces a functor $\cg \to \cg'$,
$F$ induces a functor $\Add\cg \to \Add\cg'$. Clearly, the induced
functor is essentially surjective. Let us show that it is fully
faithful. For any objects $G$ and $G'$ in $\Add\cg$ we consider the map
\[
F(G,G') : \ct(G,G') \lra \ct'(FG,FG') .
\]
By hypothesis, it is bijective if $G$ and $G'$ are in $\cg$. Let $G$
be in $\cg$ and $G' = \coprod_{i \in I}G'_i$, where $(G'_i)_{i \in I}$
is a family in $\cg$. Then, $F(G,G')$ is still bijective since we have
the following sequence of isomorphisms
\begin{eqnarray}
\nonumber \ct(G,G') & = & \ct(G,\coprod_{i \in I}G'_i) \\
\label{eqa1:i}
& = & \colim_{J \subset I}\ct(G,\coprod_{i \in J}G'_i) \\
\label{eqa1:ii}
& \iso & \colim_{J \subset I}\ct'(F(G),F(\coprod_{i \in J}G'_i)) \\
\label{eqa1:iii}
& \iso & \colim_{J \subset I}\ct'(F(G),\coprod_{i \in J}F(G'_i)) \\
\label{eqa1:iv}
& = & \ct'(F(G),\coprod_{i \in I}F(G'_i))  \\
\label{eqa1:v}
\ct'(F(G),F(G')) & \iso & \ct'(F(G),F(\coprod_{i \in I}G'_i)) \ko
\end{eqnarray}
where $J$ runs through the subsets of cardinality strictly smaller than
$\alpha$ of $I$. Here, we have used: \myref{eqa1:i} $G$ is $\alpha$-small;
\myref{eqa1:ii} $\cg$ contains $\coprod_{i \in J}G_i$ since $\cg$ is
$\alpha$-localizing; \myref{eqa1:iii} $F$ commutes with coproducts;
\myref{eqa1:iv} $F(G)$ is $\alpha$-small; \myref{eqa1:v} $F$ commutes
with coproducts. If $G'$ is in $\Add\cg$ and $G = \coprod_{i \in
  I}G_i$, where $(G_i)_{i \in I}$ is a family in $\cg$, we have
\begin{eqnarray*}
\ct(G,G') & = & \ct(\coprod_{i \in I}G_i,G')  \\
           & \iso & \prod_{i \in I}\ct(G_i,G')  \\
           & \iso & \prod_{i \in I}\ct'(F(G_i),F(G')) \\
           & \inviso & \ct'(\coprod_{i \in I}F(G_i),F(G')) \\
\ct'(F(G),F(G')) & = & \ct'(F(\coprod_{i \in I}G_i),F(G')) .
\end{eqnarray*}

\smallskip
$2^{nd}$ {\em step}: For each object $G$ in $\cg$ and each object $X$
in $\ct$, $F$ induces a bijection
\[
\ct(G,X) \lra \ct'(FG,FX) .
\]

Let $\cu$ be the full subcategory of $\ct$ formed by the objects $X$
such that $F$ induces a bijection
\[
\ct(G,X) \lra \ct'(FG,FX) \ko
\]
for each $G$ in $\cg$. Clearly, $\cu$ is a triangulated
subcategory. Let us show that $\cu$ is stable under formation of
coproducts. Let $(X_i)_{i \in I}$ be a family of objects in $\cu$. We
will show that the map
\[
\ct(G,\coprod_{i \in I}X_i) \lra \ct'(FG,F(\coprod_{i \in I}X_i) =
\ct'(FG,\coprod_{i \in I}F(X_i))
\]
is bijective. Let us show that it is surjective. Let
\[
f : FG \lra \coprod_{i \in I}F(X_i)
\]
a morphism in $\ct'$. The condition (G4) holds for the subcategory
$\cg'$ by corollary~\ref{coro:thick}. Therefore, as $F$ is an equivalence
$\cg \to \cg'$, there exists a family of objects $(G_i)_{i \in I}$ in
$\cg$ and a factorization of $f$
\[
FG \xymatrix{ \ar[rr]_g &&} \coprod_{i \in I}F(G_i) \xymatrix{
  \ar[rr]_{\coprod_{i \in I}h_i} &&} \coprod_{i \in I}F(X_i) \ko
\]
for a family of morphisms $h_i : F(G_i) \to F(X_i)$. As each $X_i$ is
in $\cu$, we have $h_i = F(k_i)$ for some morphisms $k_i : G_i \to
X_i$. Since the object
\[
\coprod_{i \in I}F(G_i) = F(\coprod_{i \in I}G_i)
\]
is in $\Add\cg'$ and $F$ induces an equivalence
\[
\Add\cg \Iso \Add\cg' \ko
\]
there exists a morphism $l : G \to \coprod_{i \in I}G_i$ such that
$F(l)$ gives $g$. Thus, $f$ is the image of the composition
\[
G \xymatrix{ \ar[rr]_l &&} \coprod_{i \in I}G_i \xymatrix{
  \ar[rr]_{\coprod_{i \in I}k_i} &&} \coprod_{i \in I}X_i
\]
under $F$. Let us show that it is injective. Let
\[
f : G \lra \coprod_{i \in I}X_i
\]
be a morphism such that $F(f) = 0$. As $\cg$ has property (G4) by
corollary~\ref{coro:thick}, we have a factorization
\[
G \xymatrix{ \ar[rr]_g &&} \coprod_{i \in I}G_i \xymatrix{
  \ar[rr]_{\coprod_{i \in I}h_i} &&} \coprod_{i \in I}X_i \ko
\]
for a family of objects $G$ in $\cg$ and a family of morphisms $h_i :
G_i \to X_i$. We have
\[
F(\coprod_{i \in I}h_i) \circ F(g) = 0 .
\]
Let us extend the morphism $\coprod_{i \in I}h_i$ and form a
distinguished triangle
\[
\coprod_{i \in I}Y_i \xymatrix{ \ar[rr]_{\coprod_{i \in I}k_i} &&}
\coprod_{i \in I}G_i \xymatrix{ \ar[rr]_{\coprod_{i \in I}h_i} &&}
\coprod_{i \in I}X_i \xymatrix{ \ar[rr] &&} \Sigma\coprod_{i \in I}Y_i
.
\]
There exists a morphism $m : FG \to \coprod_{i \in I}Y_i$ such that
\[
F(\coprod_{i \in I}k_i) \circ m = F(g) .
\]
Note that each $Y_i$ is in $\cu$ since $G_i$ and $X_i$ are in
$\cu$. By the surjectivity already shown, we have
\[
m = F(l)
\]
for a morphism $l : G \to \coprod_{i \in I}Y_i$. Thus,
\[
F(\coprod_{i \in I}k_i \circ l) = F(g) .
\]
As $G$ and $\coprod_{i \in I}G_i$ are in $\Add\cg$, it follows that
\[
\coprod_{i \in I}k_i \circ l = g .
\]
Thus,
\[
f = (\coprod_{i \in I}h_i) \circ g = (\coprod_{i \in I}h_i)
\circ (\coprod_{i \in I}k_i) \circ l = 0 .
\]

\smallskip
$3^{rd}$ {\em step}: The functor $F$ is fully faithful.

Let $Y$ be an object in $\ct$. Let $\cu$ be the full subcategory of
$\ct$ formed by the objects $X$ such that $F$ induces a bijection
\[
\ct(X,Y) \lra \ct'(FX,FY) .
\]
By the second step, $\cu$ contains $\cg$. Clearly, $\cu$ is a triangulated
subcategory. Let us show that $\cu$ is stable under formation of coproducts.
Let $(X_i)_{i \in I}$ be a family of objects in $\cu$. Then we have
\begin{eqnarray*}
\ct(\coprod_{i \in I}X_i,Y) & \iso & \prod_{i \in I}\ct(X_i,Y) \\
                      & \iso & \prod_{i \in I}\ct'(F(X_i),F(Y)) \\
                      & = & \ct'(\coprod_{i \in I}F(X_i),F(Y)) .
\end{eqnarray*}
Thus, $\coprod_{i \in I}X_i$ is indeed in $\cu$. It is easy to see
that $\cu$ contains the direct factors of its objects. So, we have
checked that $\cu$ is an $\alpha$-localizing subcategory of $\ct$ and
contains $\cg$. By proposition~\ref{prop:generators} below, the
smallest localizing subcategory containing $\cg$ is the whole category
$\ct$. It follows that $\cu = \ct$. 

\smallskip
$4^{th}$ {\em step}: The functor $F$ is essentially surjective.

The functor $F$ induces an equivalence from $\ct$ onto a localizing
subcategory $\cv$ of $\ct'$ by the third step. Indeed, $\cv$ is
triangulated, stable under coproducts and thick since $F$ is a
triangle functor commuting with coproducts. It follows that $\cv =
\ct'$, as $\cv$ contains $\cg'$, which generates $\ct'$.
\end{proof}

It remains to find conditions such that the functor $F$ of the
preceding proposition commutes with coproducts. This is made in the
following

\begin{theorem}
\label{thm:eqwgtc}
Let $\alpha$ be a regular cardinal. Let $\ct$ and $\ct'$ be two
cocomplete triangulated categories. Let $\cg \subset \ct$ and $\cg'
\subset \ct'$ be two $\alpha$-localizing subcategories, both of them
satisfying conditions (G1), (G2), (G3) for $\alpha$. Let $F : \ct \to
\ct'$ be a triangle functor. Suppose that $F$ induces a functor
\[
\cg \to \cg'
\]
which is essentially surjective and induces bijections
\[
\ct(G,X) \Iso \ct'(FG,FX)
\]
for all $G$ in $\cg$ and $X$ in $\ct$. Then $F$ is an equivalence of
triangulated categories.
\end{theorem}

\begin{remark}
We do not suppose that $F$ commutes with coproducts.
\end{remark}

\begin{proof} $1^{st}$ {\em step}: For each family $(G_i)_{i \in I}$
  in $\cg$, the morphism
\[
\coprod_{i \in I}F(G_i) \lra F(\coprod_{i \in I}G_i)
\]
is invertible.

It is sufficient to show that, for all $G'$ in $\cg'$, the map
\[
\ct'(G',\coprod_{i \in I}F(G_i) \lra \ct'(G',F(\coprod_{i \in I}G_i)
\]
is bijective since $\cg'$ verifies (G1). As $F : \cg \to \cg'$ is
essentially surjective, it is sufficient to verify this for $G' = FG$
for all $G$ in $\cg$. Let $G \in \cg$. We have
\begin{eqnarray}
\label{eqa2:1}
\ct'(FG,\coprod_{i \in I}FG_i) & = &
\colim_{J \subset I}\ct'(FG,\coprod_{i \in J}^{\ct'}FG_i)  \\
\label{eqa2:2}
& = & \colim_{J \subset I}\ct'(FG,\coprod_{i \in J}^{\cg'}FG_i) \\
\label{eqa2:3}
& = & \colim_{J \subset I}\ct(G,\coprod_{i \in J}^{\cg}G_i) \\
\label{eqa2:4}
& = & \ct(G,\coprod_{i \in I}G_i) \ko
\end{eqnarray}
where $J$ are subsets of $I$ of cardinality strictly smaller than
$\alpha$. Here, we have used: \myref{eqa2:1} $FG$ is $\alpha$-small;
\myref{eqa2:2} $\cg'$ has $\alpha$-small coproducts; \myref{eqa2:3}
$F$ induces an equivalence $\cg \to \cg'$; \myref{eqa2:4} $G$ is
$\alpha$-small. On the other hand, we have
\[
\ct'(FG,F\coprod_{i \in I}G_i) \Inviso \ct(G,\coprod_{i \in I}G_i)
\ko
\]
by the hypothesis, with $X = \coprod_{i \in I}G_i$.

\smallskip
$2^{nd}$ {\em step}: The functor $F$ induces an equivalence $\Add\cg
\to \Add\cg'$.

By the first step and the essential surjectivity of $F : \cg \to
\cg'$, $F$ induces an essentially surjective functor from $\Add\cg
\to \Add\cg'$. By hypothesis, for $G$ in $\cg$ and $X$ in $\Add\cg$,
$F$ induces a bijection
\[
\ct(G,X) \lra \ct'(FG,FX) .
\]
Let $(G_i)_{i \in I}$ be a family in $\cg$ and $X$ an object in
$\Add\cg$. Then,
\begin{eqnarray*}
\ct(\coprod_{i \in I}G_i,X) & \iso & \prod_{i \in I}\ct(G_i,X) \\
                      & \iso & \prod_{i \in I}\ct'(F(G_i),F(X)) \\
\ct'(F(\coprod_{i \in I}G_i),F(X)) & = & \ct'(\coprod_{i \in
  I}F(G_i),F(X)) .
\end{eqnarray*}
Thus, $F$ restricted to $\Add\cg$ is fully faithful.

\smallskip
$3^{rd}$ {\em step}: The functor $F$ commutes with coproducts.

Let us consider the diagram
\[
\xymatrix{
\ct \ar[r]^F \ar[d]_<>(0.5){h_\ct} & \ct' \ar[d]^<>(0.5){h_\ct'} \\
\coh(\Add\cg) \ar[r]^\sim_{F^*} & \coh(\Add\cg') .
}
\]
We will show that it commutes up to isomorphism. Let $X$ be an object
in $\ct$ and let
\[
h(G_1) \lra h(G_0) \lra h(X) \lra 0
\]
be a projective presentation, where $G_1$, $G_0$ are in
$\Add\cg$. Then, for all $G$ in $\cg$, we obtain an exact sequence
\[
\ct(G,G_1) \lra \ct(G,G_0) \lra \ct(G,X) \lra 0 .
\]
Therefore, the sequence
\[
\ct'(FG,FG_1) \lra \ct'(FG,FG_0) \lra \ct'(FG,FX) \lra 0
\]
is exact (since isomorphic to the first). It follows that the sequence
\[
h(FG_1) \lra h(FG_0) \lra h(FX) \lra 0
\]
is exact (since the objects $h(FG)$, $G \in \cg$, form a family of
projective generators of $\coh(\Add\cg')$). Thus,
\[
F^*(h(X)) = \cok(h(FG_1) \lra h(FG_0))
\]
is indeed canonically isomorphic to $h(FX)$. To conclude, note that
$F^*$ (which is an equivalence!) and $h_\ct$ commute with coproducts
and that $h_\ct'$ detects the isomorphisms.

\smallskip
$4^{th}$ {\em step}: The claim follows thanks to the preceding
proposition~\ref{prop:eqwgtc}.
\end{proof}

\section{Thick subcategories and localization of triangulated categories}
\label{s:tsubcatloc}

We recall now some known results about the localizations of
triangulated categories and about their thick subcategories, before
stating the most important theorem of this section concerning the
localization of {\em well generated} triangulated categories. For
complete proofs of the cited results, we refer to Neeman's book
\cite[Ch.~2, p.~73]{Neeman99} and the classical \cite[Ch.~2.2, p.~111-133]{Verdier96}.

\subsection{Localization of triangulated categories}
\label{ss:loctc}
We begin with a collection of properties of the {\em triangle quotient}
\cite[Ch.~2.2, p.~111-133]{Verdier96} of triangulated categories.

\begin{proposition}
\label{prop:locgeneral}
Let $\ct$ be a triangulated category with arbitrary coproducts, let
$\Phi$ be a {\em set} of morphisms in $\ct$ and $\cn$ the smallest
triangulated subcategory of $\ct$ containing the $\cone (s)$, with
$s \in \Phi$, stable under arbitrary coproducts. Then the following
assertions hold:
\begin{itemize}
\item[a)] $\ct/\cn$ is a triangulated category and admits arbitrary
  coproducts;
\item[b)] the canonical functor $Q : \ct \to \ct/\cn$ commutes with
  all coproducts;
\item[c)] the morphisms $Q(s)$ are invertible for all $s\in\Phi$;
\item[d)] if $F : \ct \to \cs$ is a triangle functor, where $\cs$
  is a triangulated category which admits all coproducts and the
  functor $F$ commutes with all coproducts and makes every $s \in
  \Phi$ invertible in $\cs$, then $F=\ol F \circ Q$ for a unique
  coproduct preserving triangle functor $\ol F : \ct/\cn \to \cs$;
\item[e)] more precisely, if $\cs$ is a triangulated category with
  arbitrary coproducts, there is an isomorphism of categories
\[
\xymatrix{
\cF un_{cont}(\ct/\cn,\cs) \ar[r]^\sim & \cF un_{cont,\Phi}(\ct,\cs)
\ko
}
\]
where $\cF un_{cont}$ is the category of triangulated functors
commuting with arbitrary coproducts, and $\cF un_{cont,\Phi}$ is the
category of the functors in $\cF un_{cont}$ which have the additional
property of making all $s \in \Phi$ invertible.
\end{itemize}
\end{proposition}

\begin{proof}
See Chapter 2 in \cite{Neeman99} or \cite[Ch.~2.2, p.~111-133]{Verdier96}. We give only an
argument for the commutativity of $Q$ with coproducts because it is a
general one, useful in other contexts. Let $\prod_I\ct$ be the product
category of copies of $\ct$ indexed by $I$.
%\ie the category with objects the families $(X_i)_{i\in I}$ of
%objects of $\ct$ and as morphisms the obvious ones induced
%componentwise by those in $\ct$. Consider the diagonal functor
%$\Delta : \ct \to \prod_I\ct$. $X \mapsto (X)_I$ and $(f : X \to Y)
%\mapsto ((f)_I : (X)_I \to (Y)_I)$.
Using the universality of coproducts it is easy to check that the
functor $\coprod_{i\in I} : \prod_I\ct \to \ct$ which takes a family
$(X_i)_{i \in I}$ to the coproduct $\coprod_{i \in I}X_i$ is left
adjoint to the diagonal functor $\Delta$. It is clear that
$\Delta(\Phi) \subseteq \prod_I(\Phi)$ and that $\coprod_{i\in I}
(\prod_I(\Phi)) \subseteq \Phi$. Therefore, the pair
$\coprod_{i\in I} \dashv \Delta$ induces the following commutative
diagram
\[
\xymatrix{
\prod_I\ct \ar[rrr]^<>(0.5){can} \ar@<-0.5ex>[d]_{\coprod_{i\in I}}
&&& (\prod_I\ct)[(\prod_I\Phi)^{-1}] \ar@2{-}[r] &
\prod_I(\ct[\Phi^{-1}]) \ar@<-0.5ex>[d]_{\coprod_{i\in I}}   \\
\ct \ar[rrrr]^<>(0.5){can} \ar@<-0.5ex>[u]_{\Delta} &&&&
\ct[\Phi^{-1}] \ar@<-0.5ex>[u]_{\Delta}
}
\]
which entails the required commutativity of $Q$ with all the coproducts
over the set $I$. Of course, this construction is possible for every set
$I$.
\end{proof}

The functor $Q$ is usually called {\em (canonical) quotient functor}, even if $\ct$ does not have coproducts. In general, the morphisms between two objects in a triangle quotient do not form a set. However, this is the case if the quotient functor $Q$ admits a right adjoint $Q_\rho$, because $Q_\rho$ is automatically fully faithful.
%In this article, we will mainly deal with a particular kind of localization,
%called Bousfield localization \cite[ch. 8]{Neeman99}.
% 
% \begin{definition}
% \label{def:Bouloc}
% Let $\ct$ be a triangulated category with small Hom-sets and $\cn$ a
% triangulated subcategory of $\ct$. If the canonical quotient functor
% $Q$ admits a right adjoint $\Qr$, we say that $\Qr$ is the {\em
%   Bousfield localization functor} for the pair $\cn \subseteq
% \ct$.
% \end{definition}
%
%If we do not require any condition on the $\Hom$-sets, we have the
%following, more general, definition.

\begin{definition}
\label{def:loc}
Let $\ct$ and $\ct'$ be triangulated categories. A triangle functor $F
: \ct' \to \ct$ is a {\em localization functor} if it is fully faithful and
admits a left adjoint functor.
\end{definition}
%The two definitions are equivalent if in the last, we require that the
% category $\ct$ has small Hom-sets. In the following we will usually
% only write {\em localization} and {\em localization functor},
% specifying if the category has coproducts everytime.

If $F : \ct' \to \ct$ is a localization functor and $F_\lambda$ is left adjoint, then $F_\lambda$ induces an equivalence from the triangle quotient $\ct \slash \ker(F_\lambda)$ to $\ct'$. Via this equivalence, $F_\lambda$ identifies with the quotient functor $\ct \to \ct \slash \ker(F_\lambda)$ and $F$ with its right adjoint.

\subsection{Some thick subcategories of triangulated categories}
\label{ss:subcat}
Let us recall that a full triangulated subcategory of a triangulated
category is called {\em thick}, ({\em {\'e}paisse} in the French literature, {\em satur{\'e}e} in the original definition in Verdier's thesis \cite[2.2.6, p.~114]{Verdier96}) if it contains the direct factors of its objects. We remark that this property is automatically verified if the triangulated category has countable coproducts, since in this case idempotents splits (see \cite{BalmerSchlichting01}, \cite[Prop.~1.6.8, p.~65]{Neeman99} for definitions and properties of idempotents in triangulated categories). Now we give
definitions and notations about some important subcategories of a
triangulated category $\ct$ with arbitrary coproducts and suspension functor
$\Sigma$. The best reference for this material is \cite[Ch.~3-4]{Neeman99}. We
recall that a full triangulated subcategory $\cs$ of $\ct$ is called
{\em localizing} if it is closed under arbitrary coproducts. It is called
$\alpha${\em -localizing}, for a given regular cardinal $\alpha$, if it
is thick and closed under $\alpha${\em -coproducts} of its objects,
\ie coproducts of objects formed over a set of cardinality strictly
smaller than $\alpha$. We write $\Sa$ for the smallest
$\alpha$-localizing subcategory of $\ct$ containing $\cs$, where $\cs$
is a set or a class of objects in $\ct$ and $\alpha$ a {\em regular}
cardinal. Note that in the above definitions the requirement that the
subcategories are thick is necessary only for the case $\alpha =
\aleph_0$, since for $\alpha > \aleph_0$ these subcategories are
automatically thick as we underlined at the beginning of the
section. In his book \cite[Ch.~3-4]{Neeman99} Neeman shows the very important
properties of the previous subcategories and of the subcategories of the
$\alpha$-small objects $\ct^{(\alpha)}$ and that of the
$\alpha$-compact objects $\ta$. They are triangulated,
$\alpha$-localizing and thick subcategories of $\ct$ for
$\alpha>\aleph_0$. There is the following filtration: if $\alpha \leq
\beta$ then $\ta \subseteq \tb$. If $\cs$ contains a good set of
generators for $\ct$, then $\ct = \bigcup_\alpha \!\Sa$, where
$\alpha$ runs over all regular cardinals. If $\ct$ is a well generated
triangulated category, then $\ct = \bigcup_\alpha \ta$, where $\alpha$
runs over all regular cardinals.

We have the

\begin{theorem}[{\cite[Lemma 3.2.10, p.~107]{Neeman99}}] Let $\beta$ be an infinite cardinal. Let $\ct$ be a triangulated category closed under the formation of coproducts of fewer than $\beta$ of its objects. Let $\cn$ be a $\beta$-localising subcategory of $\ct$. Then $\ct \slash \cn$ is
closed with respect to the formation of coproducts of fewer than $\beta$ of its
objects, and the universal functor $F : \ct \to \ct \slash \cn$ preserves coproducts.
\end{theorem}

and its

\begin{corollary}[{\cite[Cor. 3.2.11, p.~110]{Neeman99}}] If $\ct$ is a
  triangulated category with all coproducts and $\cn$ is a localizing subcategory of $\ct$, then $\ct/\cn$ is a triangulated category which admits all coproducts and the universal functor $\ct \to \ct/\cn$ preserves coproducts.
\end{corollary}

Now we state one of the major results in the theory of triangulated
categories. This result has a long story (see for example \cite{Brown62}, \cite{Neeman98}, \cite[Ch.~10]{KashiwaraSchapira05}), which comes from algebraic topology. We state it in the modern and general form Neeman gives it in his book, see \cite[Thm.~1.17, p.~16]{Neeman99} and \cite[Thm.~8.3.3, p.~282]{Neeman99} for a more general statement and the proof.

\begin{theorem}
\label{thm:brown}
(Brown representability). Let $\ct$ be a well-generated triangulated
category. Let $H$ be a contravariant functor $H : \ct\op \to \ca b
$. The functor $H$ is representable if and only if it is
cohomological and takes coproducts in $\ct$ to products of abelian
groups.
\end{theorem}

Let us now clarify the meaning of two similar but different
notions. Let $\ct$ be a triangulated category and $\cg$ a set of
objects in $\ct$. We say that $\ct$ is {\em generated by $\cg$} or,
equivalently, that {\em $\cg$ generates $\ct$} if $\ct = \Ggen$. In
contradistinction, we say that $\cg$ is a {\em generating set} for $\ct$
if condition (G1) of definition~\ref{def:wgtcKrause} holds for the subset
$\cg$ of $\ct$. Let us recall condition (G1): an object $X \in \ct$ is
zero if and only if $\ct(G,X) = 0$ for all $G$ in $\cg_0$ (we always
assume that $\Sigma\cg_0 = \cg_0$). The former notion is stronger than
the latter. If $\cg$ generates $\ct$, then $\cg$ is a generating set
for $\ct$, whereas the converse is not true in general but holds if $\cg$ is
assumed further to be a {\em $\aleph_1$-perfect generating set for $\ct$} in the sense of Neeman (\cf \cite[Ch.~8, Def.~8.1.2, p.~273]{Neeman99}). Moreover, we give another link between the two notions, which is useful because it covers the case of well generated triangulated categories.

\begin{proposition}
\label{prop:generators}
Let $\ct$ be a well generated triangulated category and $\cg$ a
generating set for $\ct$, \ie condition {\em (G1)} holds for
$\cg$. Then $\cg$ generates $\ct$, \ie $\Ggen = \ct$.
\end{proposition}

\begin{proof}
Let us call $\cn$ the subcategory $\Ggen$. Since $\cn$ is a localizing
subcategory generated by the {\em set} $\cg$, it is well generated by
corollary~\ref{coro:locwgtc} below. Then, the Brown representability
theorem \myref{thm:brown} holds for $\cn$. Therefore, for each object
$X \in \ct$, the functor $\Hom_\ct(-,X)|_\cg : \Ggen\op \to \ca b$,
which is cohomological and sends coproducts into products, is
representable. For each object $X$ in $\ct$ there exists an object
$X_\cn$ in $\cn$ such that $\Hom_\ct(-,X)|_\cg \iso
\Hom_\cn(-,X_\cn)$. Thus, we have obtained a functor $i_\rho$ right
adjoint to the fully faithful inclusion: $i : \cn \to \ct$. Consider
now, for every $X \in \ct$, the distinguished triangle in $\ct$
\[
ii_\rho X \lra X \lra Y \lra \Sigma ii_\rho X .
\]
Applying to the triangle the covariant functor $\Hom_\ct(iN,-)$, $N$
an object of $\cn$, we obtain a long exact sequence of abelian
groups. Consider the part corresponding to the input triangle: The map
from the first term $\Hom_\ct(iN,ii_\rho X)$ to the second term
$\Hom_\ct(iN,X)$ is easily seen to be an isomorphism, since $i$ is fully
faithful and $i_\rho$ is its right adjoint. Similarly, the map from the
fourth to the fifth term is an isomorphism. Therefore, the third group
$\Hom_\ct(iN,Y)$ must be zero for all $N \in \cn$. This forces the object
$Y$ to lie in $\cn^\perp$. But $\cn^\perp$ is zero. Indeed, condition (G1)
holds for $\cg$, \ie $\cg^\perp = 0$. Thus, the inclusion $\cg \subseteq
\cn$ gives $\cn^\perp \subseteq \cg^\perp = 0$, \ie $\cn^\perp =
0$. Therefore, we have $Y = 0$. By the triangle above, this means
$ii_\rho X \iso X$, for all $X$ in $\ct$. It follows that $i$ is an
equivalence of categories, which gives $\ct = \Ggen$.
\end{proof}
%\begin{remark}
%\label{rmk:Tria}
%We recall for the convenience of the reader that Bernhard Keller
%provides in \cite{Keller06} a characterization of the subcategory
%$\Ggen$ in the particular case where $\ct$ is the derived category of
%a small DG category $\ca$. In this case he defines the {\em
%  triangulated category $\tria(\ca)$ associated with $\ca$} to be the
%closure in $\cd(\ca)$ of the set of representable functors $\Ac$, $A
%\in \ca$, under shifts in both directions and extensions, and the {\em
%  category of perfect objects $\per(\ca)$} to be the closure of
%$\tria(\ca)$ under passage to direct factors in $\cd(\ca)$. He shows
%that the category $\per(\ca)$ equals the category of compact objects
%$(\DA)^c$. We can continue and define the {\em triangulated category
%  $\Tria(\ca)$ associated with $\ca$} to be the closure of $\per(\ca)$
%under small coproducts and the {\em triangulated category
%  $\Tria_\alpha(\ca)$ associated with $\ca$}, where $\alpha$ is a
%regular cardinal, to be the closure of $\per(\ca)$ under small
%$\alpha$-coproducts. The class formed by the objects of $\per(\ca)$ is
%a set, since $\ca$ is small, and is a generating set for $\DA$. Then,
%$\Tria(\ca)$ equals $\DA$ by
%proposition~\ref{prop:generators}. Moreover, $\Tria_\alpha(\ca)$
%clearly equals $\DaA$.
%\end{remark}

\subsection{Localization of well generated triangulated categories}
\label{ss:locwgtc}
In this section, we will state a Theorem about particular localizations
of well generated triangulated categories, those which are triangle
quotients by a subcategory {\em generated by a set}. One could obtain
this result using Thomason's powerful Theorem \cite[Key Proposition 5.2.2, p.~338]{ThomasonTrobaugh90} in its generalized form given by Neeman in \cite[Thm.~4.4.9, p.~143]{Neeman99}. Neeman himself does this in \cite{Neeman01} in proving that the derived category of a Grothendieck category is always a well generated triangulated category. We will give a more detailed and slightly different proof in order to make clear the machinery behind Thomason-Neeman's Theorem. Before doing this task we recall the key ingredient of the proof.

\begin{theorem}[{\cite[Ch.~4, Thm.~4.3.3, p.~131]{Neeman99}}]
\label{thm:factorization}
Let $\ct$ be a triangulated category with small coproducts. Let $\beta$
be a regular cardinal. Let $\cs$ be some class of objects in $\tb$. Let
$X$ be a $\beta$-compact object of $\ct$, \ie $X\in\tb$, and let $Z$ be
an object of $\Sgen$. Suppose that $f : X \to Z$ is a morphism in $\ct$.
Then there exists an object $Y \in \Sb$ so that $f$ factors as $X \to Y
\to Z$.
\end{theorem}

\begin{proof}
As in Neeman's book, since that proof uses only the facts that $\tb$ is
a $\beta$-localizing triangulated subcategory of $\ct$, that all the
objects in $\tb$ are $\beta$-small, and that condition (G4) holds for
$\tb$. These properties are also valid for Krause's definition of
$\tb$.
\end{proof}

The power of this property is seen at once, since it is the key to
obtain the following results.

\begin{corollary}[{\cite[Ch.~4, Lemma 4.4.5, p.~140 for item a) and Lemma 4.4.8, p.~142 for item b)]{Neeman99}}]
\label{coro:corofactorization}
Let $\ct$ be a triangulated category with small coproducts. Let $\cs$
be some class of objects in $\ta$ for some infinite cardinal
$\alpha$. Let $\beta \geq \alpha$ be a regular cardinal. Then:
\begin{itemize}
\item[a)] if $\Sgen=\ct$, then the inclusion $\Sb\subseteq\tb$ is an
  equality;
\item[b)] let $\cn=\Sgen$. Then there is an inclusion $\cn \cap \tb
\subseteq \nb$.
\end{itemize}
\end{corollary}

\begin{proof}
a) Let $X$ be an object of $\tb$ and consider the identity map $\id_X
: X \to X$. As X is at the same time in $\tb$ and in $\Sgen$ we can
apply the theorem~\ref{thm:factorization} and factor $\id_X$ through
some object $Y\in\Sb$. Thus the object $X$ is a direct factor of $Y$.
Since $\Sb$ is thick, we have $X \in \Sb$.

b) Let $K$ be an object of $\cn \cap \tb$. Then $K$ is $\beta$-small
as an object of $\cn$ since the inclusion $\cn \subseteq \ct$ commutes
with coproducts. Now, let $K \to \coprod_{i\in I} X_i$ be a morphism,
where the objects $X_i$ belong to $\cn$. It factors through a morphism
$\coprod_{i\in I} f_i : \coprod_{i\in I} K_i \to \coprod_{i\in I}
X_i$, where the objects $K_i$ belong to $\tb$. By the theorem above,
each morphism $K_i \to X_i$ factors through an object $K'_i$ belonging
to $\Sb \subseteq \cn\cap\tb$. Therefore the class $\cn \cap \tb$
satisfies (G4) in $\cn$ and we obtain the required inclusion.
\end{proof}

The next proposition states some useful properties of the images in
the quotient category $\ct/\cn$ of the maps of the subcategories $\Gb$
of $\ct$ under the canonical quotient functor $Q$.

\begin{proposition}
\label{prop:Q-image}
Let $\alpha$ be a regular cardinal. Let $\ct$ be a triangulated
category with small coproducts, generated by a class of objects $\cg
\subseteq \ta$. Let $\cs$ be an arbitrary class of objects in $\ta$
and $Q$ the canonical quotient functor
\[
Q : \ct \to \ct/\Sgen .
\]
Let $\beta \geq \alpha$ be a regular cardinal. Then:
\begin{itemize}
\item[a)] each morphism $u : Q(G) \to Q(X)$, where $G$ is an object of
      $\Gb$ and $Q(X)$ an arbitrary object of $\ct/\cn$, is the
      equivalence class of a diagram in $\ct$  \\
      \[
      \xymatrix@ur{
      G & G' \ar[l]^\sim \ar[d]  \\
      & X \ko
      }
      \]
      where the object $G'$ belongs to $\tb = \Gb$ and the arrow
      $\iso$ means a morphism whose image under $Q$ is invertible; in
      particular, the morphisms from $Q(G)$ to $Q(X)$ in $\ct \smash \Sgen$
      form a set if $\cg$ is a set;
\item[b)] the image of $\Gb$ under the (restriction of the) functor
      $Q$ is a full triangulated subcategory of $\ct/\cn$;
\item[c)] if $\beta$ is uncountable, then $\QGb$ equals $Q(\Gb)$. If
      $\beta$ is countable, then $\QGb$ equals the closure of $Q(\Gb)$
      under taking direct factors.
\end{itemize}
\end{proposition}

\begin{proof}
a) Let $u : Q(G) \to Q(X)$ be a morphism in $\ct/\cn$. It is the
equivalence class of a `roof' diagram in $\ct$  \\
\[
\xymatrix@ur{
G & T \ar[l]^\sim \ar[d]  \\
& X \ko
}
\]
where the object $T$ belongs to $\ct$. We can form the distinguished
triangle
\[
\xymatrix{
N & G \ar[l] & T \ar[l]_\sim & \Sigma^{-1}N \ar[l] \ko
}
\]
where $N$ and $\Sigma^{-1}N$ lie in $\Sgen$. The object $G$ is
$\beta$-compact in $\ct$. Therefore, we can apply
theorem~\ref{thm:factorization} to the morphism $N \la G$ and factor
it as
\[
N \lla N' \lla G \ko
\]
where $N'$ belongs to $\Sb$. The class $\Sb$ is contained in $\tb$,
since $\cs$ is contained in $\tb$ by the hypothesis. Therefore, we can
complete the morphism $N' \la G$ to a distinguished triangle in $\tb$
\[
\xymatrix{
N' & G \ar[l] & G' \ar[l]_\sim & \Sigma^{-1}N' \ar[l]
}
\]
and deduce a map of distinguished triangles \\
\[
\xymatrix@ur{
N' \ar[d] & G \ar[l] \ar@{=}[d] & G' \ar[l]^\sim \ar@{-->}[d]^\sim
\ar@/^1.3pc/@{.>}[dd] & \Sigma^{-1}N' \ar[l] \ar[d]   \\
N. & G \ar[l] \ar@{~>}[rd] & T \ar[l]^\sim \ar[d] & \Sigma^{-1}N
\ar[l]  \\
&& X
}
\]
adding the morphism $G' \to T$. The wavy arrow stands for the given
morphism $Q(G) \to Q(X)$ in $\ct/\cn$, whereas the dotted arrow is the
composition $G' \to T \iso X$. The roof diagrams $G \inviso G' \to X$
and $G \inviso T \to X$ are clearly equivalent. We have supposed that
$\ct$ has small coproducts and that $\Ggen = \ct$, with $\cg$
contained in $\ta$, hence in $\tb$. Therefore, $\Gb = \tb$ by point a)
of corollary~\ref{coro:corofactorization}. This shows that $G'$ also
lies in $\Gb$.

\smallskip

b) Clearly, the image of $\Gb$ under $Q$ is stable under $\Sigma$ and $\Sigma^{-1}$.
We have to show that it is stable under forming cones. Let $G_1$ and $G_2$ be two
objects of $\Gb$ and $u$ a morphism from $QG_1$ to $QG_2$. By part a), the morphism
$u$ equals the equivalence class of a diagram
\[
\xymatrix@ur{
G_1 & G_1' \ar[l]^\sim \ar[d]^v  \\
& G_2 \ko
}
\]
where $G_1'$ belongs to $\Gb$. Therefore, the cone $C$ on $v$ still belongs to $\Gb$.
Clearly, the cone on $u$ is isomorphic to $Q(C)$, which still belongs to the image
under $Q$ of $\Gb$.

\smallskip
%c) Suppose that the regular cardinal $\beta$ is strictly greater than
%$\aleph_0$. The inclusion $Q(\Gb) \subseteq \QGb$ follows from the
%fact that $Q$ is a triangle functor and commutes with arbitrary
%coproducts. For the reverse inclusion, we use the fact that the
%triangulated subcategory $Q(\Gb)$ of $\ct/\cn$ has
%$\alpha$-coproducts, hence countable coproducts. Therefore, it contains
%the direct factors of its objects (\cf~\ref{ss:subcat}). Moreover,
%$Q(\Gb)$ contains the image of $\cg$ under $Q$, hence is
%$\beta$-localizing by point b). The claim for $\beta > \aleph_0$
%follows by the definition of $\QGb$, the smallest $\beta$-localizing
%triangulated subcategory of $\ct/\cn$ containing $Q\cg$. If $\beta =
%\aleph_0$, the subcategory $Q(\langle\cg\rangle_{\aleph_0})$ does not
%contain the direct factors of its objects. Thus, $\langle Q\cg
%\rangle_{\aleph_0}$ equals the closure of $Q(\langle \cg
%\rangle_{\aleph_0})$ under taking direct factors. In this case,
%$Q(\langle\cg\rangle_{\aleph_0})$ is usually said {\em dense} in
%$\langle Q\cg \rangle_{\aleph_0}$ (\cf \cite[1.4, 1.5]{Thomason97}).

c) Let $\cu$ be the closure of $Q(\Gb)$ under taking direct factors.
We claim that $\cu$ equals $\QGb$ for all $\beta \geq \alpha$. Indeed,
we have $Q(\Gb) \subseteq \QGb$ since $Q$ is a triangle functor and
commutes with arbitrary coproducts. It follows that $\cu \subseteq
\QGb$ since $\QGb$ is thick. For the reverse inclusion, we notice that
$\cu$ contains $Q\cg$, that it is a triangulated subcategory since
$Q(\Gb)$ is a triangulated subcategory (by b), and that it is thick (by
definition). We have thus proved the claim for countable $\beta$. Now
suppose $\beta$ is uncountable. Then $Q(\Gb)$ is a triangulated
subcategory stable under forming countable coproducts. Therefore, it
is stable under taking direct factors (\cf~\ref{ss:subcat}) and thus
equals $\cu = \QGb$.
\end{proof}

Now we can state the most important theorem of this section. This
theorem has been inspired by Neeman's generalization to well generated
categories \cite[Thm.~4.4.9, p.~143]{Neeman99} of Thomason-Trobaugh's theorem
\cite[Key Proposition 5.2.2, p.~338]{ThomasonTrobaugh90}.

\begin{theorem}
\label{thm:locwgtc}
Let $\ct$ be an $\alpha$-compactly generated triangulated category and
$\cg$ a set of good generators for $\ct$, contained in $\ta$. Let $\cs$
be a {\em set} of objects contained in $\tg$, for some fixed regular
cardinal $\gamma$. Let $\cn = \Sgen$ and $Q$ the canonical quotient
functor
\[
Q : \ct \to \ct/\cn .
\]
\begin{itemize}
\item[a)] The localizing triangulated subcategory $\cn$ is the union
      \[
      \cn = \bigcup_{\delta\geq\gamma} \nd \ko
      \]
      where $\delta$ runs through the regular cardinals. Equivalently, $\cn$ is
      given by the same union as above, formed over all regular cardinals;
\item[b)] the subcategory $\cn$ is $\delta$-compactly generated for
      all regular cardinals $\delta \geq \gamma$ by the set $\Sg$;
\item[c)] the subcategory $Q(\Gb)$ equals $\QGb$ for $\beta >
      \aleph_0$ and its closure under taking direct factors equals
      $\QGb$ for $\beta = \aleph_0$;
\item[d)] the quotient category $\ct/\cn$ is a $\delta$-compactly
      generated triangulated category for all regular cardinals
      $\delta \geq \beta$, where $\beta = \sup(\alpha,\gamma)$, with
      set of good generators $Q(\Gb)$.
\end{itemize}
\end{theorem}

\begin{proof}
It is clearly sufficient to prove b) for $\delta = \gamma$ and c) for
$\delta = \beta$.

\smallskip
a) The triangulated category $\ct$ is well generated. Therefore, it
is the union over all the regular cardinals $\sigma$ of its
subcategories $\ts$ \cite[Corollary of Thm. A]{Krause01}. We know
from the hypothesis that $\cs \subseteq \tg$, hence $\cs \subseteq
\cn \cap \tg$. Clearly, $\Sg \subseteq \cn \cap \tg$, since $\Sg$ is
the smallest $\gamma$-localizing subcategory of $\ct$ containing the
set $\cs$. Moreover, $\cn \cap \tg \subseteq \ngamma$ by point b) of
corollary~\ref{coro:corofactorization}. Thus, we have the following
sequence of inclusions:
\[
\cs \subseteq \Sg \subseteq \cn \cap \tg \subseteq \ngamma .
\]
Therefore, for each regular cardinal $\delta \geq \gamma$, we obtain
$\Sd = \cn \cap \td = \nd$, by point a) of
corollary~\ref{coro:corofactorization}. The claim now follows by the
equalities
\[
\cn = \cn \cap \ct = \cn \cap(\bigcup_\lambda \tl) =
\bigcup_\lambda(\cn \cap \tl) = \bigcup_{\delta\geq\gamma} \nd =
\bigcup_\lambda \nl .
\]
The two last equalities hold since the set of
the subcategories $\nl$ is filtered over regular cardinals. This means
that $\na \subseteq \nb$ if $\alpha \leq \beta$, for all regular
cardinals $\alpha$ and $\beta$.

\smallskip
b) The isomorphism classes of the objects of the subcategory $\Sg$
form a set, since it is explicitly constructed from the objects in
$\cs$, which is also a set. Moreover, $\Sg$ is stable under shifts
because it is triangulated. Let us show condition (G1). Let $Y$ be
an object of $\cn$ such that $\Hom_\cn(X,Y)=0$ for all $X$ in $\Sg$.
Then, it is easy to check that this equality holds for $X$ in
$\Sgen$. In particular it holds for $X=Y$. Hence $Y$ vanishes.
Therefore, condition (G1) holds for $\Sg$. By the proof of point a), $\Sg =
\ngamma$. Therefore, conditions (G2) and (G3) trivially hold by the
definition of $\ngamma$.

\smallskip
%MISTAKE: the Hom could be a class between two objects! d) Let us
%show that the hypotheses of proposition~\ref{prop:Q-image} hold. We
%begin by showing that the category $\ct$ has small $\Hom$-sets. We
%have supposed that the triangulated category $\ct$ is
%$\alpha$-compactly generated and that the set $\cg$ is contained in
%$\ta$. Therefore, $\Gd$ equals $\td$ for all regular cardinals $\delta
%\geq \alpha$, by the point a) of corollary~\ref{coro:corofactorization}.
%Therefore, we can write
%\[
%\ct = \bigcup_{\delta\geq\alpha} \td = \bigcup_{\delta\geq\alpha} \Gd
%\ko
%\]
%where $\delta$ runs over regular cardinals. The subcategories $\Gd$
%are essentially small because they are explicitly constructible from
%the set $\cg$. Therefore, each pair $X$, $Y$ of objects in $\ct$ must
%lie in $\langle\cg\rangle_{\ol\delta}$ for some $\ol\delta \geq
%\alpha$. Thus, we can write
%\[
%\Hom_\ct(X,Y) = \Hom_{\langle\cg\rangle_{\ol\delta}}(X,Y) \ko
%\]
%since $\langle\cg\rangle_{\ol\delta}$ is a full subcategory of
%$\ct$. This shows that the category $\ct$ has small
%$\Hom$-sets (NO!!).
c) All the conditions of proposition~\ref{prop:Q-image} hold. Thus,
this point results from part c) of proposition~\ref{prop:Q-image}.

\smallskip
d) The subcategory $\cn$ is well generated by point b). Thus, the Brown
representability theorem \myref{thm:brown} holds for $\cn$ and we conclude that the inclusion $i$ of $\cn$ into $\ct$ admits a right adjoint $i_\rho$ as in the proof of proposition \ref{prop:generators}. Now this implies that the quotient functor $Q : \ct \to \ct \slash \cn$ admits a right adjoint $Q_\rho$ (which takes an object $X$ to the cone of the adjunction morphism $ii_\rho X \to X$). The functor $Q_\rho$ is a localization functor \myref{def:loc}. Thus, it is fully faithful. Let us sum up the situation in the following diagram,
\[
\xymatrix{
\cn \ar@{ (->}@<+0.2ex>[rr]^<>(0.5)i && \ct \ar@<+0.4ex>[rr]^<>(0.5)Q \ar@<+0.4ex>[ll]^<>(0.5){i_\rho} && \ct/\cn \ar@{ (->}@<+0.2ex>[ll]^<>(0.5)\Qr .
}
\]
We have to show that the conditions (G1), (G2) and (G3) of
definition~\ref{def:wgtcKrause} hold for $Q(\Gb)$. We begin by
observing that the sets $\cg$ and $\cs$ are both contained in $\tb$,
since we have chosen $\beta = \sup(\alpha,\gamma)$. The condition
(G1) holds even for the smaller set $Q\cg$, hence for $Q(\Gb)$.
Indeed, suppose $\Hom_{\ct/\cn}(Q\cg,X) = 0$, for an arbitrary object
$X$ in $\ct/\cn$. By the adjunction, this is equivalent to
$\Hom_\ct(\cg,\Qr(X)) = 0$. The condition (G1) holds for the set $\cg$
in $\ct$ and implies $\Qr(X) = 0$. Therefore, $X = Q\Qr X = 0$, since
$Q\Qr$ is naturally equivalent to the identity endofunctor of
$\ct/\cn$. Thus, condition (G1) holds for the set $Q\cg$. The
subcategory $Q(\Gb)$ contains its $\beta$-coproducts because $Q$
commutes with all coproducts and its objects form a set. Therefore,
conditions (G2) and (G4) are equivalent for $Q(\Gb)$ (\cf \cite[Lemma
4]{Krause01}). Let us now simultaneously show that conditions (G4) and
(G3) hold for $Q(\Gb)$. Consider a morphism $u : Q(G) \to
\coprod_{i\in I} X_i$, where $G$ is an arbitrary object in $\Gb$. We
know from point a) of proposition~\ref{prop:Q-image} that $u$ is the
equivalence class of a diagram in $\ct$  \\
\[
\xymatrix@ur{
G & G' \ar[l]^\sim \ar[d]^f  \\
& \coprod_{i\in I} X_i \ko
}
\]
where the object $G'$ belongs to $\tb = \Gb$. The conditions (G3) and
(G4) also hold for $\Gb$, by corollary~\ref{coro:thickwgtc}. Therefore,
there exists a set $J \subset I$ of cardinality strictly smaller than
$\beta$ and a set of morphisms $(f_i : G_i \to X_i)_{i\in I}$, where
$G_i$ lies in $\Gb$ for all $i \in J$, so that the morphism $f : G'
\to \coprod_{i\in I} X_i$ factors through $\coprod_{i\in J}X_i$ (G3)
\[
\xymatrix@ur{
G' \ar[ddrr]^<>(0.5)f \ar[dd]  \\\\
\coprod_{i\in J} X_i \ar[rr] && \coprod_{i\in I} X_i
}
\]
and through the morphism $\coprod_{i\in I}f_i$ (G4)
\[
\xymatrix@ur{
G' \ar[ddrr]^<>(0.5)f \ar[dd]  \\\\
\coprod_{i\in I} G_i. \ar[rr]_{\coprod_{i\in I}f_i} && \coprod_{i \in
  I} X_i
}
\]
The image under $Q$ of the the last two diagrams shows that the
morphism $u$ factors in $\ct/\cn$ in the same way. Therefore,
conditions (G3) and (G4) hold for $Q(\Gb)$.
\end{proof}

\begin{remark}
\label{rmk:bigbeta}
The construction of the cardinal $\beta$ in the preceding proof is
not optimized at all. In spite of the constructive proof, this
result will be useful mainly for existence problems.
\end{remark}

The next corollary is a result about the localization of well
generated categories obtained by inverting a {\em set} of arrows,
implicitly contained in Neeman's book \cite{Neeman99}.

\begin{corollary}
\label{coro:locwgtc}
Let $\ct$ be a well generated triangulated category and $\cn$ a
localizing triangulated subcategory of $\ct$ generated by a {\em set}
of objects $\cs$. Then $\cn$ and $\ct/\cn$ are well generated
triangulated categories.
\end{corollary}

\begin{proof}
Take the coproduct of all the objects in $\cs$. Since $\cs$ is a set,
the coproduct will be in $\tg$ for some regular cardinal $\gamma$.
Therefore, we have $\cs \subseteq \tg$, because $\tg$ is thick in
$\ct$ and so contains the direct factors of its objects. Now apply
theorem \ref{thm:locwgtc}.
\end{proof}

\section{The $\alpha$-continuous derived category}
\label{s:acdc}

In this section, we construct the {\em $\alpha$-continuous} derived category of
a homotopically $\alpha$-cocomplete DG category. This construction enjoys a
useful and beautiful property. Given a homotopically $\alpha$-cocomplete (\cf
below) pretriangulated DG category $\ca$, we will show that its
$\alpha$-continuous derived category $\DaA$ is $\alpha$-compactly generated by the free DG modules. The categories $\DaA$ will be the prototypes of the $\alpha$-compactly generated algebraic DG categories. We use the notations of \cite{Keller06}.

\begin{definition}
Let $\alpha$ be a regular cardinal and $\ca$ a small DG
$k$-category. We assume that $\ca$ is {\em homotopically
  $\alpha$-cocomplete}, \ie that the category $H^0(\ca)$ admits all
$\alpha$-small coproducts. For each $\alpha$-small family $(A_i)_{i
  \in I}$ of objects of $\ca$, we write
\[
\coprodh_{i \in I}A_i
\]
for their coproduct in $H^0(\ca)$. Each DG functor $M : \ca^{op} \to
\cc_{dg}(k)$ induces a functor $H^0 M : (H^0(\ca))^{op} \to \ch(k)$
and so we have a canonical morphism
\[
(H^* M)(\coprodh_{i \in I}A_i) \xymatrix{ \ar[r] &} \prod_{i \in
  I}(H^* M)(A_i).
\]
Let $\DA$ be the derived category \cite{Keller94} of $\ca$. It is a
triangulated category (\cf \cite{Keller94}, \cite{Keller06}). The
{\em $\alpha$-continuous derived category} $\DaA$ is defined as the
full subcategory of $\DA$ whose objects are the DG functors $M$ such
that, for each $\alpha$-small family of objects $(A_i)_{i \in I}$ of
$\ca$, the canonical morphism above is invertible.
\end{definition}

\begin{remark}
All the small $k$-linear DG categories which are {\em
  $\alpha$-cocomplete}, \ie admit all $\alpha$-small coproducts, are
homotopically $\alpha$-cocomplete. A partial converse is given in remark~\ref{rmk:thesis} below.
\end{remark}
%\begin{remark}
%\label{rmk:set}
%The DG category $\ca$ is supposed small, \ie its objects form a set.
%Therefore, the objects of the subcategory $\cg$ of representables of
%$\DA$, form a set. Hence the subcategories $(\DA)^\alpha$ are
%essentially small for all regular cardinals $\alpha$, since they
%equal the subcategories $\Ga$ starting from an $\alpha$ great enough
%\myref{coro:corofactorization}. Indeed, these last are essentially
%small since their objects are explicitly constructed from the
%objects of $\cg$. In particular, these subcategories contain all the
%$\alpha$-small coproducts of representables, whose isomorphism
%classes therefore form a set.
%\end{remark}

This definition describes $\DaA$ as a subcategory of $\DA$. One can
give an equivalent definition in terms of a localization of $\DA$,
which yields a category $\DAN$ triangle equivalent to $\DaA$. For
this, let us define some {\em sets} of morphisms in the category of
DG modules $\CA$ (\cf \cite{Keller94}, \cite{Keller06}). We recall
that the notation $A^\wedge$ means $\Hom_\ca(-,A)$. Let $\Sigma_0$
be the set of all morphisms of $\CA$
\[
\sigma_\lambda : \coprod_{i \in I}\Aic \xymatrix{ \ar[r] &}
(\coprodh_{i \in I}A_i)^\wedge \ko
\]
%\[
%\sigma_\lambda : \coprod_{i \in I}(\Aic) \lra (\coprodh_{i \in
%  I}A_i)^\wedge \ko
%\]
where $\lambda$ ranges over the {\em set} $\Lambda$ of all
families $(A_i)_{i \in I}$ in $\ca$ of cardinality strictly smaller
than $\alpha$. We define $\Sigma$ to be the set of cofibrations
\cite{Quillen67}, \cite{Hovey99}, \cite{Hirschhorn03} between
cofibrant DG modules \cite{Keller06}
\[
\begin{bmatrix}
\sigma_\lambda \\ inc
\end{bmatrix} :
\coprod_{i \in I}\Aic \xymatrix{ \ar@{ >->}[r] &} (\coprodh_{i \in
  I}A_i)^\wedge \oplus I(\coprod_{i \in I}\Aic) \ko
\]
where $\lambda \in \Lambda$ and, for each object $X$, the morphism
$inc : X \xymatrix{ \ar@{ >->}[r] &} IX$ is the inclusion of $X$ into
the cone over its identity morphism. We can also consider the set
\[
\cm = \{ \Nl \xymatrix{ \ar@{ >->}[r] &}
I\Nl~|~\Nl=\cone(\sigma_\lambda) \vir \lambda \in \Lambda \}.
\]
The cones over the morphisms in $\Sigma_0$, $\Sigma$ or $\cm$ generate
the same localizing subcategory $\cn$ of $\DA$ because the objects
$IX$ are contractible and thus become zero objects in $\DA$. The
quotient functor
\[
\DA \lra \DAN
\]
induces an equivalence
\[
\DaA \Iso \DAN.
\]
Following \cite{Drinfeld04}, we say that a DG category $\ca$ is {\em
  pretriangulated} if the essential image of the Yoneda functor is a
triangulated subcategory of the derived category $\DA$. In the case of
pretriangulated DG categories, the definition of {\em
  quasi-equivalence of DG categories} of \cite{Keller06} specializes to
the following.

\begin{definition}
Let $\ca$ and $\ca'$ be pretriangulated DG categories. A {\em DG
  functor} $F : \ca \to \ca'$ is a {\em quasi-equivalence of
  pretriangulated DG categories} if the induced triangle functor
$H^0(F) : H^0(\ca) \to H^0(\ca')$ is an equivalence of triangulated
categories.
\end{definition}
%The next proposition states that a homotopically $\alpha$-cocomplete
%pretriangulated DG category contains $\alpha$-small coproducts up to a
%quasi-equivalence.
%We suspect that the answer to the following question is negative.

\begin{remark}
\label{rmk:thesis}
We will show in \cite{PortaThesis} that if $\ca$ is a homotopically $\alpha$-cocomplete pretriangulated $\DG$ category, then there exists a quasi-equivalence $\ca \to \ca'$, where $\ca'$ is a pretriangulated DG category which is {\em $\alpha$-cocomplete}. This establishes the link between this article and \cite{Tabuada07}.
\end{remark}

We now come to the result which motivated the definition of the
$\alpha$-continuous derived category.

\begin{theorem}
\label{thm:rightimage}
Let $\ca$ be a homotopically $\alpha$-cocomplete pretriangulated
$\DG$ category. The $\alpha$-continuous derived category of $\ca$ is
$\alpha$-compactly generated by the images of the free $\DG$ modules
$\Ac$, $A \in \ca$. More precisely, the full subcategory $\cg$ of
$\DaA$ formed by the images of the free DG modules $\Ac$, $A \in \ca$,
is a triangulated subcategory satisfying conditions (G1), (G2) and
(G3) of definition~\ref {def:wgtcKrause}.
\end{theorem}

\begin{remark}
\label{rmk:aleph_0}
We prove the theorem in the case where $\alpha$ is strictly greater
than $\aleph_0$, the case $\alpha = \aleph_0$ being trivial. In fact,
$\aleph_0$-coproducts are finite coproducts. Thus, the morphisms
$\sigma_\lambda$ above are isomorphisms already in $\DA$, and
$\cd_{\aleph_0}\ca$ equals $\DA$.
\end{remark}

\begin{proof} This proof depends heavily on
  theorem~\ref{thm:locwgtc}. Therefore, let us explain how the
  notations correspond. The triangulated category
  $\ct$ is $\DA$. The {\em set} $\cs$ is formed by the cones on the
  following morphisms
\[
\sigma_\lambda : \coprod_{i \in I}(\Aic) \xymatrix{ \ar[r] &}
(\coprodh_{i \in I}A_i)^\wedge \ko
\]
where $\lambda$ ranges over the {\em set} $\Lambda$ of all families
$(A_i)_{i \in I}$ in $\ca$ of cardinality strictly smaller then
$\alpha$. The set $\cg$ is formed by the free DG modules $\Ac$, $A
\in \ca$. It is contained in $\ct^{\aleph_0}$, whereas $\cs$ is
contained in $\ta$. We have $\beta = \sup(\aleph_0,\alpha) =
\alpha$. Let $\cn$ be $\Sgen$ and $Q$ the projection functor
\[
Q : \DA \lra \DaA \eqiso \ct/\cn .
\]
Then, according to theorem~\ref{thm:locwgtc}, the $\alpha$-continuous
derived category $\DaA$ is $\alpha$-compactly generated by
$\QGa$. Hence, the claim of the theorem is equivalent to the
following claim: $\QGa = Q\cg$ and the functor $Q$ induces an
equivalence $\cg \iso Q\cg$. We begin with the equivalence $\cg \iso
Q\cg$. It amounts to the same as to show that the functor $Q|_\cg$
is fully faithful. We know from the proof of the point d) of
theorem~\ref{thm:locwgtc} that $Q$ admits a right adjoint $\Qr$
\myref{def:loc}. From the general theory of Bousfield localizations
\cite[Ch.~9]{Neeman99}, we have that $Q|_{\cn^\perp} : \cn^\perp \to \DaA$ is an
equivalence of triangulated categories. In particular $Q|_{\cn^\perp}$ is fully
faithful. Therefore, it is sufficient to show that $\cg$ is
contained in $\cn^\perp$. By definition, $\cn$ is the localizing
subcategory generated by the cones $\cone(\sigma_\lambda)$, which we
call $C_{(A_I)}$. We have to show that each $\Ac \in \cg$ is right
orthogonal to the objects $C_{(A_I)}$. By applying the cohomological
functor $\Hom_\DA(-,\Sigma^n\Ac)$, $n \in \Z$, to the distinguished
triangle
\[
\coprod_{i \in I}(\Aic) \lra (\coprodh_{i \in I}A_i)^\wedge \lra
C_{(A_I)} \lra \Sigma\coprod_{i \in I}(\Aic) \ko
\]
it is clear that it is sufficient to show that the natural morphism
\[
\Hom_\DA(\coprod_{i \in I}(\Aic),\Sigma^n\Ac) \lla
\Hom_\DA((\coprodh_{i \in I}A_i)^\wedge,\Sigma^n\Ac)
\]
is an isomorphism for all $n \in \Z$. This follows from the following
sequence of isomorphisms
%\begin{align*}
%\Hom_\DA(\coprod_{i \in I}(\Aic),\Sigma^n\Ac) & \iso
%\prod_{i \in I}\Hom_\DA(\Aic,\Sigma^n\Ac) \\
%& = \prod_{i \in I}H^n\Hom_\ca(A_i,A) \\
%& \inviso H^n(\prod_{i \in I}\Hom_\ca(A_i,A)) \\
%\Hom_\DA((\coprodh_{i \in I}A_i)^\wedge,\Sigma^n\Ac)
%& \iso H^n(\Hom_\ca(\coprodh_{i \in I}A_i,A)).
%\end{align*}

\begin{eqnarray*}
\Hom_\DA(\coprod_{i \in I}(\Aic),\Sigma^n\Ac) & \iso &
\prod_{i \in I}\Hom_\DA(\Aic,\Sigma^n\Ac) \\\\
& = & \prod_{i \in I}H^n\Hom_\ca(A_i,A) \\\\
& \inviso & H^n(\prod_{i \in I}\Hom_\ca(A_i,A)) \\
\Hom_\DA((\coprodh_{i \in I}A_i)^\wedge,\Sigma^n\Ac)
& \iso & H^n(\Hom_\ca(\coprodh_{i \in I}A_i,A)).
\end{eqnarray*}
The second and the last isomorphisms are justified by the following
one
\[
\Hom_\DA(\Ac,\Sigma^n B^\wedge) = H^n\Hom_\ca(A,B) \ko
\]
valid for all $A$ and $B$ in $\ca$.
%which is isomorphism~\ref{eq:isorepder} of section~\ref{ss:dercat}
%in the case $M = \Ac$.
The third isomorphism is the fact that cohomology commutes
with formation of products. For the fourth, we observe that the
natural homomorphism
\[
\Hom_\ca(\coprodh_{i \in I}A_i,A) \lra \prod_{i \in I}\Hom_\ca(A_i,A)
\]
is a homotopy equivalence, by the definition of $\coprodh$. Therefore,
it becomes invertible in cohomology.

It is trivial that $Q\cg$ is stable under shifts. Moreover, it is
automatically thick for $\alpha > \aleph_0$ \myref{ss:subcat}. To
prove  that $Q\cg$ equals $\QGa$ it is then sufficient to show that
$Q\cg$ is closed under $\alpha$-coproducts and extensions. We have
\[
\coprod_{i \in I}(Q\Aic) \Iso Q(\coprod_{i \in I}\Aic) \Iso
Q((\coprodh_{i \in I}A_i)^\wedge) \ko
\]
where the cardinality of $I$ is strictly smaller than $\alpha$ and the
last isomorphism holds by the construction of $\cn$. This shows that
$Q\cg$ is closed under formation of $\alpha$-coproducts. Finally,
$Q\cg$ is stable under extensions in $\DaA$. Indeed, $\cg$ is stable
under extensions in $\DA$ and hence in $\cn^\perp$, since we have
shown that $\cg$ is contained in $\cn^\perp$. We have also seen that
the restriction $Q|_{\cn^\perp}$ is an equivalence of the categories
$\cn^\perp$ and $\DaA$. It follows that $Q\cg$ is stable under
extensions in $\DaA$.
\end{proof}

\section{The Popescu-Gabriel theorem for triangulated categories}
\label{s:PGforTC}

\subsection{Algebraic triangulated categories}
\label{ss:algtriacat}
Let us recall that an exact category $\ce$ \cite{Keller90},
\cite{Quillen73} is a Frobenius category \cite{Happel88} if it has
enough injectives, enough projectives, and the two classes of the
injectives and projectives coincide. For all pairs of objects $X$,
$Y$ of $\ce$, let $I_\ce(X,Y)$ be the subgroup of the abelian group
$\Hom_\ce(X,Y)$ formed by the morphisms which factor over an
injective-projective object of $\ce$. The {\em stable category of
$\ce$} \cite{Happel88}, written $\ul \ce$, is the category which
has the same objects as $\ce$ and the morphisms
\[
\Hom_{\ul \ce}(X,Y) = \Hom_\ce(X,Y) / I_\ce(X,Y).
\]

\begin{definition} \cite{Keller06}
An {\em algebraic} triangulated category is a $k$-linear triangulated
category which is triangle equivalent to the stable category $\ul \ce$
of some $k$-linear Frobenius category $\ce$.
\end{definition}

The class of algebraic triangulated categories is stable under taking
triangulated subcategories and forming triangulated localizations (up
to a set-theoretic problem). Examples abound since categories of
complexes up to homotopy are algebraic. Therefore, the categories
arising in homological contexts in algebra and geometry are
algebraic. The area where one often encounters non algebraic
triangulated categories is topology. In particular the stable homotopy
category of spectra is not algebraic. More examples can be found in
section 3.6 of \cite{Keller06}.

\subsection{The main theorem}
\label{ss:mainthm}
We recall \cite{Keller94} that a {\em graded category} over a
commutative ring $k$ is a $k$-linear category $\cb$ whose morphism
spaces are $\Z$-graded $k$-modules
\[
\Hom_\cb(X,Y) = \coprod_{p \in \Z} \Hom_\cb(X,Y)^p
\]
such that the composition maps
\[
\Hom_\cb(X,Y) \ten_k \Hom_\cb(Y,Z) \lra
\Hom_\cb(Y,Z)
\]
are homogeneous of degree $0$, for all $X$, $Y$, $Z$ in $\cb$. Now we
can state and prove the main theorem of this paper.

\begin{theorem}
\label{thm:mainthm}
Let $\ct$ be a triangulated category. Then the following statements
are equivalent:
\begin{itemize}
\item[($i$)] $\ct$ is algebraic and well generated;
\item[($ii$)] there is a small $\DG$ category $\ca$ such that $\ct$ is
      triangle equivalent to a localization of $\DA$ with respect to a
      localizing subcategory generated by a {\em set} of objects.
\end{itemize}
Moreover, if $\ct$ is algebraic and $\alpha$-compactly generated, and
$\cu \subset \ct$ is a full triangulated subcategory stable under
$\alpha$-small coproducts and such that conditions (G1), (G2) and
(G3) of definition~\ref{def:wgtcKrause} hold for $\cu$, then there is
an associated localization functor \myref{def:loc} $\ct \to \DA$ for
some small $\DG$ category $\ca$ such that $H^*(\ca)$ is equivalent to
the graded category $\cu_{gr}$ whose objects are those of $\cu$ and
whose morphisms are given by
\[
\cu_{gr}(U_1,U_2) = \bigoplus_{n\in\Z} \ct(U_1,\Sigma^n U_2) .
\]
\end{theorem}

\begin{proof}
$(ii) \Longrightarrow (i)$ : $\ct$ is a localization of $\DA$, \ie
there is a fully faithful functor
\[
\xymatrix{
\ct \ar@{ (->}@<-0.05ex>[rr]^<>(0.5)\Ft && \DA \ko
}
\]
admitting a left adjoint functor. The category $\DA$ is
algebraic. Triangulated subcategories of algebraic categories are
algebraic, implying that $\ct$ is algebraic, too. Moreover, $\DA$ is
compactly generated by the {\em set} ($\ca$ is small)
\[
\{~X^\wedge[n]~|~n \in \Z \vir X \in \ca~\}
\]
thanks to the isomorphism
\[
\Hom_\DA(X^\wedge[n],M) \iso H^{-n}(M(X)) \ko
\]
where $M$ is a DG module and $X$ is an object of $\ca$ (\cf
\cite{Keller94}, \cite{Keller06}). Therefore, $\ct$ is well generated
by corollary~\ref{coro:locwgtc}, since it is assumed to be a
localization generated by a set of the $\aleph_0$-compactly generated
category $\DA$.

\smallskip
$(i) \Longrightarrow (ii)$ : for the sake of clarity, we will give the
proof of this implication in several steps, after making the main
construction.

\smallskip
Let $\ct$ be an algebraic, well generated triangulated category, \ie
$\ct$ is equivalent to $\ul \ce$ for some Frobenius category $\ce$. By
the definition of well generated triangulated category (in the sense
of Krause), there are a regular cardinal $\alpha$ and a set of
$\alpha$-good generators $\cg_0 \subseteq \ct$ such that $\Sigma\cg_0
= \cg_0$ and the conditions (G1), (G2) and (G3) of
definition~\ref{def:wgtcKrause} hold. Let $\cg$ be the closure of the
set $\cg_0$ under extensions and $\alpha$-coproducts. The set $\cg$ is
stable under the suspension functor $\Sigma$ of $\ct$ and under its
inverse. Therefore, it is a small triangulated subcategory of
$\ct$. Let us recall and summarize the properties which hold for
$\cg$.
\begin{itemize}
\item[(G0)] The set $\cg$ is a small full triangulated subcategory of
      $\ct$, stable under the formation of all $\alpha$-small
      coproducts;
\item[(G1)] the set $\cg$ is a generating set for $\ct$: An object
  $X \in \ct$ is zero if $\Hom_\ct(G,X)=0$ for all $G$ in $\cg$;
\item[(G3)] all the objects $G \in \cg$ are {\em $\alpha$-small}: For
  each family of objects $X_i$, $i \in I$, of $\ct$, we have
  $\Hom_\ct(G,\coprod_I X_i) = \colim_{J \subset I}
  \Hom_\ct(G,\coprod_J X_i)$, where the sets $J$ have cardinality
  strictly smaller than $\alpha$;
\item[(G4)] for each family of objects $X_i$, $i\in I$, of $\ct$, and
  each object $G \in \cg$, each morphism
\[
G \lra \coprod_{i\in I} X_i
\]
factors through a morphism $\coprod_{i\in I} \phi_i$: $\coprod_{i \in
  I} G_i \to \coprod_{i \in I} X_i$, with $G_i$ in $\cg$ for all $i
  \in I$.
\end{itemize}
Condition (G0) clearly holds for $\cg$. Condition (G3) of
definition~\ref{def:wgtcKrause} has just been rewritten using
colimits. Condition (G4) holds for $\cg$ by
proposition~\ref{prop:dense}. Note that conditions (G2) and (G4) are
equivalent for $\cg$. Indeed, we can apply \cite[Lemma 4]{Krause01},
since the set $\cg$ has $\alpha$-coproducts and its objects are
$\alpha$-small. 

We may assume that the category $\ce$ is of the form $Z^0(\tilde\ce)$
for an exact DG category $\tilde\ce$ by the argument of the proof of
theorem 4.4 of \cite{Keller94}. Let us recall that a DG category $\ca$
is an {\em exact $\DG$ category} \cite{Keller99} if the full
subcategory $Z^0(\ca)$ of $\CA$ formed by the image of the Yoneda
functor is closed under shifts and extensions (in the sense of the
exact structure of \cite{Keller94}). Then, $H^0(\ca)$ becomes a
triangulated subcategory of $\ch(\ca)$ and the subcategory of the
representable functors becomes a triangulated subcategory of $\DA$.
Thus, an exact DG category is also a pretriangulated DG category
(\cf section~\ref{s:acdc}). Let us now define a small full DG
subcategory $\ca \subset \tilde\ce$ as follows. For each isomorphism
class of objects of $\cg$, we choose a representative $G$ and we
denote by $A_G$ the same object considered in the category
$\tilde\ce$. By definition, these objects $A_G$ are objects of $\ca$.
Then, clearly, the category $H^0(\ca)$ is a full subcategory of
$H^0(\tilde\ce) = \ul\ce = \ct$ and it is equivalent to $\cg$ by the
functor sending $A_G$ to $G$. In particular, $H^0(\ca)$ is a
triangulated category and it admits all $\alpha$-small coproducts.
Thus, $\ca$ is a homotopically $\alpha$-cocomplete pretriangulated
DG category. We define the functor
\[
F : \ct \lra \DaA
\]
by sending an object $X$ of $\ct = H^0(\tilde\ce)$ to the DG module
$FX$ taking $A_G \in \ca$ to $\Hom_{\tilde\ce}(G,X)$. A priori, $FX$
lies in $\DA$. Let us show that it belongs in fact to the full
subcategory $\DaA$ of $\DA$. Let $A_{G_i}$, $i \in I$, be an
$\alpha$-small family in $\ca$. Then the coproduct $\coprod^{H^0}_{i
  \in I}A_{G_i}$ of the $A_{G_i}$ in $H^0(\ca)$ is isomorphic to
$A_{\coprod_{i \in I}G_i}$. Thus, we have a quasi-isomorphism
\[
(FX)(\coprod^{H^0}_{i \in I}A_{G_i}) = \Hom_{\tilde\ce}(\coprod^\ct_{i
  \in I}G_i,X) \lra \prod_{i \in I}\Hom_{\tilde\ce}(G_i,X) = \prod_{i
  \in I}(FX)(A_{G_i}) \ko
\]
induced by liftings to $\ce = Z^0(\tilde\ce)$ of the canonical
morphisms $G_j \to \coprod^\ct_{i \in I}G_i$ in $\ct$, respectively by
representatives in $Z^0(\ca)$ of the canonical morphisms $A_{G_j} \to
\coprod^{H^0}_{i \in I}A_{G_i}$ in $H^0(\ca)$. For $A_G \in \ca$, we
have
\[
FG = \Hom_{\tilde\ce}(-,G) = \Hom_\ca(-,A_G) = A_G^\wedge \ko
\]
which shows that $F$ induces an essentially surjective functor from
$\cg$ to the full subcategory of the $A_G^\wedge$ in $\DaA$. For $A_G$
in $\ca$ and $X$ in $\ct$, we have
\begin{eqnarray*}
\Hom_\DaA(FG,FX) & = & \Hom_\DaA(A_G^\wedge,FX)  \\
& = & \Hom_\DA(A_G^\wedge,FX)  \\
& = & H^0(FX(A_G))  \\
\Hom_\ct(G,X) & = & H^0(\tilde\ce(G,X)).
\end{eqnarray*}
We would like to apply theorem~\ref{thm:eqwgtc} to conclude that $F$
is a triangle equivalence: In the notations of
theorem~\ref{thm:eqwgtc}, we take $\ct = \ct$, $\cg = \cg$, $\ct' =
\DaA$ and $\cg'$ to be the full subcategory on the objects
$A_G^\wedge$ in $\DaA$. By theorem~\ref{thm:rightimage}, $\ct'$ and
$\cg'$ do satisfy the hypothesis of theorem~\ref{thm:eqwgtc} and so
$F$ is indeed a triangle equivalence.

\smallskip
Now suppose that $\ct$ is an algebraic well generated triangulated
category. Let $\cu \subset \ct$ be a full small subcategory as in the
last assertion of the theorem. Then the conditions (G0), (G1) and (G3)
above hold for $\cg = \cu$. Moreover, condition (G4) holds for $\cg =
\cu$ by \cite[Lemma 4]{Krause01}. Therefore, we can construct a DG
category $\ca$ and an equivalence $F : \ct \iso \DaA$ as above in the
proof of the implication from i) to ii). Moreover, $H^*(\ca)$ equals
$\cu_{gr}$. Indeed, both have the same objects and we have
\begin{eqnarray*}
H^n(\ca)(A_{G_1},A_{G_2}) & = &
\ch\ca(A_{G_1}^\wedge,\Sigma^n(A_{G_2}^\wedge))  \\
& = & \ch\ca(A_{G_1}^\wedge,(\Sigma^n A_{G_2})^\wedge) \\
& = & H^0(\ca)(A_{G_1},\Sigma^n A_{G_2}) \\
\cu_{gr}(G_1,G_2)^n  & = & \cu(G_1,\Sigma^n G_2) .
\end{eqnarray*}
\end{proof}
%We give a summarizing commutative diagram in order to have a
%global vision of the situation described in detail in the last proof:
%\[
%\xymatrix{
%\ct \ar@/^/[drr]^\Ft \ar[dr]^F \ar@/_/[ddr]_{h_\ct} \\
%& \DaA \ar@{ (->}@<+0.3ex>[r]^\Qr \ar[d]^{H^0_\alpha} & \DA
%\ar@<+0.3ex>[l]^Q \ar[d]^{H^0} \\
%& \Mod_\alpha \cg \ar@{ (->}@<+0.3ex>[r]^R & \Mod\cg
%\ar@{.>}@<+0.3ex>[l]^L
%}
%\]
%where $H^0_\alpha$ is the restriction of $H_0$ to $\DaA$, $\Mod\cg =
%\{ M : \cg \op \longrightarrow \ca b \}$ and $\Mod_\alpha \cg = \{
%M~|~M \mbox{ takes $\alpha$-small coproducts to products} \}$. Note
%that the functor $L$ is not exact in general and then $\Mod_\alpha
%\cg$ is not a localization of $\Mod\cg$.
%
%\smallskip
%% Corollaries: G-P for triang. cat.s; case aleph_0 gives old thm.
%The next corollary will give a triangulated compact version of the
%famous Gabriel-Popescu theorem \cite{PopescuGabriel64} as we stated it
%in the introduction \myref{thm:PG}. Note the structural equivalence
%with the abelian case. Here is the hint to say that Keller's algebraic
%categories are the right generalisation of Grothendieck's categories
%(without generators) and Neeman's well generated triangulated
%categories are the right ones for additive categories with
%generators.
%
%\begin{corollary}
%Let $\ct$ be a well generated algebraic triangulated category. Then
%the following statements are equivalent:
%\begin{itemize}
%\item[($i$)] $G\in\ct$ is an arbitrary generator of $\ct$;
%\item[($ii$)] $\RHom(G,-) : \ct \longrightarrow \cd(A) \vir A =
%\RHom(G,G)$ is a localization.
%\end{itemize}
%\end{corollary}

If $\ct$ is compactly generated we recover a result obtained by
B. Keller in \cite[Thm.~4.3]{Keller94}:

\begin{corollary}
Let $\ct$ be an algebraic triangulated category. Then the following
statements are equivalent:
\begin{itemize}
\item[($i$)] $\ct$ is compactly generated;
\item[($ii$)] $\ct$ is equivalent to the derived category $\DA$ for
some small $\DG$ category $\ca$.
\end{itemize}
\end{corollary}

\begin{proof}
See remark~\ref{rmk:aleph_0}
\end{proof}

\subsection{Application}
\label{ss:applex}
We apply theorem \ref{thm:mainthm} to a certain class of
subcategories of algebraic triangulated categories we are going to
define.

\begin{definition}
Let $\ct$ be an algebraic triangulated category which is triangle
equivalent to the stable category of the Frobenius category $\ce$ and
admits arbitrary coproducts. Let $\tilde\ce$ be a DG category (not
necessarily small) such that $H^0(\tilde\ce)$ is triangle equivalent
to $\ct$. Given a subcategory $\cg$ of $\ct$, let $\tilde\cg$ be the DG
subcategory of $\ca$ with the same objects as $\cg$. Thus, the category
$H^0(\tilde\cg)$ is isomorphic to $\cg$. We will say that $\cg$ is a
{\em compactifying} subcategory of $\ct$ if it is small and the functor
\[
\ct \longrightarrow \cd\tilde\cg \vir X \longmapsto
\Hom_{\tilde\ce}(-,X)|_{\tilde\cg}
\]
is fully faithful.
\end{definition}

For example, W. T. Lowen and M. Van den Bergh proved in
\cite[Ch.~5]{LowenVandenBergh04} that, given a Grothendieck category
$\ca$ with a generator $G$, the one-object subcategory $\cg = \{ G \}$
of the derived category $\DA$ is a compactifying subcategory. For this
reason we call such a generator $G$ {\em compactifying}.

\begin{theorem}
Let $\ct$ be a well generated algebraic triangulated category. Then
there is a regular cardinal $\alpha$ such that the subcategory
$\sk(\tb)$ formed by a system of representatives of the isomorphism
classes of $\tb$ is compactifying for each regular cardinal $\beta
\geq \alpha$.
\end{theorem}

\begin{proof}
Suppose that $\alpha$ is the first regular cardinal such that $\ct =
\Tagen$. This cardinal exists because the category $\ct$ is well
generated. For each $\beta \geq \alpha$, the subcategory $\sk(\tb)$
is small and satisfies conditions (G1), (G2) and (G3) of
definition~\ref{def:wgtcKrause} by definition of the subcategory $\tb$
and the filtration by increasing regular cardinals. Now the claim
follows from the last part of theorem~\ref{thm:mainthm}.
\end{proof}

%\bibliographystyle{amsplain}
%\bibliography{stanKeller}

\def\cprime{$'$}
\providecommand{\bysame}{\leavevmode\hbox to3em{\hrulefill}\thinspace}
\providecommand{\MR}{\relax\ifhmode\unskip\space\fi MR }
\providecommand{\MRhref}[2]{%
  \href{http://www.ams.org/mathscinet-getitem?mr=#1}{#2}
}
\providecommand{\href}[2]{#2}

\end{document}